\newif\ifproofs
\proofsfalse

\documentclass[graybox]{svmult}

\usepackage{fancyvrb}
\usepackage{graphicx,psfrag,caption,subcaption,tikz}
\usepackage[cmex10]{amsmath}
\usepackage{hyperref}
\usepackage{booktabs}
\usepackage{color}
\usepackage{psfrag}
\usepackage[toc,page]{appendix}
\definecolor{price}{rgb}{0.000, 0.500, 0.000}
\definecolor{flow}{rgb}{0.500, 0.500, 0.000}
\definecolor{pmnt}{rgb}{0.000, 0.000, 0.500}

\newcommand{\BEAS}{\begin{eqnarray*}}
\newcommand{\EEAS}{\end{eqnarray*}}
\newcommand{\BEQ}{\begin{equation}}
\newcommand{\EEQ}{\end{equation}}
\newcommand{\BIT}{\begin{itemize}}
\newcommand{\EIT}{\end{itemize}}

\newcommand{\eg}{{\it e.g.}}
\newcommand{\ie}{{\it i.e.}}

\newcommand{\reals}{{\mbox{\bf R}}}

\newcommand{\pmin}{p_{\rm min}}
\newcommand{\pmax}{p_{\rm max}}
\newcommand{\pfix}{p_{\rm fix}}

\newcounter{oursection}

\usetikzlibrary{calc,patterns,decorations.pathmorphing,decorations.markings,arrows}

\tikzset{
  position/.style args={#1:#2 from #3}{
    at=(#3.#1), anchor=#1+180, shift=(#1:#2)
  }
}

\title{Dynamic Energy Management}

\author{Nicholas Moehle, Enzo Busseti, Stephen Boyd, and Matt Wytock}
\institute{Nicholas Moehle \at Department of Mechanical Engineering, Stanford University, \email{nicholasmoehle@gmail.com}
\and Enzo Busseti \at Department of Management Science and Engineering, Stanford University, \email{ebusseti@stanford.edu}
\and Stephen Boyd \at Department of Electrical Engineering, Stanford University, \email{boyd@stanford.edu}
\and Matt Wytock \at Gridmatic, Inc., \email{matt@gridmatic.com}
}

\date{\today}

\begin{document}
\maketitle

\tikzstyle{device}=[draw,thick,fill=black!10,node distance=3cm,
                    minimum height=1cm,minimum width=2cm]
\tikzstyle{net}=[draw,circle,thick,node distance=3cm,
                 minimum width=.7cm, inner sep=0cm]
\tikzstyle{classicterm} = [draw, very thick]
\tikzstyle{terminal} = [draw, very thick, -latex']
\tikzstyle{termnoarrow} = [draw, very thick]
\tikzstyle{bus} = [draw, line width=3pt]

\newcommand{\abstracttext}{
We present a unified method, based on convex optimization,
for managing the power produced and consumed by a network of devices over time. 
We start with the simple setting of optimizing power flows in
a static network, and then proceed to the case of optimizing 
dynamic power flows, \ie, power flows that change with time over a horizon.
We leverage this to develop a real-time control strategy,
model predictive control, 
which at each time step solves
a dynamic power flow optimization problem, using 
forecasts of future quantities such as demands, capacities, or prices,
to choose the current power flow values.
Finally, we consider a useful extension of model predictive control
that explicitly accounts for uncertainty in the forecasts.
We mirror our framework with an object-oriented software
implementation, an open-source Python library for planning and controlling
power flows at any scale.
We demonstrate our method with various examples.
Appendices give more detail about the package, and describe some
basic but very effective methods for constructing forecasts from 
historical data.
}

\abstract*{\abstracttext}
\abstract{\abstracttext}





\section{Introduction}
We present a general method for planning power production, consumption, conversion, 
and transmission throughout a network of interconnected devices.
Our method is based on convex optimization.  It provides power flows
that meet all the device constraints, as well as conservation of power between
devices, and minimizes a total cost associated with the devices.
As a by-product, the method determines the locational marginal price for power at
each point on the network where power is exchanged.

In the simplest setting we ignore time and consider static networks.
In the next simplest setting, we optimize power flows for multiple time periods,
over a finite time horizon, which allows us to include ramp rate constraints, energy
storage devices, and deferrable loads.
We leverage this to develop a real-time control method, 
model predictive control, which uses forecasts of unknown quantities and optimization
over a horizon to create a plan, the first step of which is used or executed in
each time period.
It is well known that, despite uncertainty in the forecasts, 
model predictive control often works reasonably well.
Finally, we consider an extension of model predictive control that
explicitly handles uncertainty in the forecasts by considering several 
possible scenarios, and creating a full contingency plan for each one,
coupled by the requirement that the power flows in the first period of 
each contingency must be the same.


In addition to providing optimal power flows,
our method computes the locational marginal price of power over the network.
These prices can be used as the basis of a system of payments,
in which each device is paid for the power it produces,
or pays for the power it consumes, and each transmission line
or conversion device is paid its service.
We show that, under this payment scheme,
the optimal power flows maximize each individual device's profit, \ie, 
the income from payments to it minus the cost of operating the device.
This means that the optimal power flows are not only socially optimal,
but also provide an economic equilibrium,
\ie, there is no incentive for any device to deviate from the optimal power flows
(in the absence of price manipulation).

Our exposition is accompanied by \texttt{cvxpower},
a Python software implementation of our method,
which is available at
\url{http://github.com/cvxgrp/cvxpower}.
\texttt{cvxpower} is an object-oriented package
which provides a declarative language
for describing and optimizing power networks.
Object-oriented software design is well suited for building
complex applications with many inter-operating components,
whose users need not understand the internal details of these 
components.
In this sense, we aim to abstract away the technical details of
the individual devices in the network,
as well as the underlying optimization problem,
allowing users to focus on modeling.
More advanced users can extend our software framework,
for example, by defining and implementing a new device.

Most of the ideas, and much of the material in this paper has appeared in 
other works or is well known.  Our contribution is to assemble it all into one
coherent framework, with uniform notation and an organization of ideas that shows
how a very basic method of optimizing static power flows generalizes
naturally to far more complex settings.
We also note that more sophisticated work has appeared on closely related topics, 
such as using convex optimization relaxations to solve (nonconvex) AC power flow
optimization problems, or advanced forms of robust model predictive control.
Some of this work is discussed in the section below on related work,
as well as in the main body of this paper.

\subsection{Related work}
Mathematical optimization has been used to manage electric power grids
for nearly a century.
Modern overviews of the field are provided by
\cite{wood2012power} and \cite{taylor2015convex},
and many examples by \cite{baldick2006applied}.
For an overview of economic and financial issues related to energy markets
see \cite{harris2011electricity}.

\paragraph{\emph{Optimal dispatch.}}
One of the earliest applications of optimization to power systems
is the optimal dispatch problem,
which considers the problem of planning the operation of multiple
generators in order to meet power demand.
This method dates back to 1922.
(See \cite{davison1922dividing}.)
and a good classical treatment can be found in
\cite{steinberg1943economy}.
This early work was based on the \emph{incremental rate method},
which involves solving the optimality conditions
of a convex optimization problem by hand,
using graphical methods.
(These conditions are (\ref{e-static-opt-cond}) in our formulation,
for a problem with a single net.)
For a typical, modern formulation,
see \cite[Ch.~3]{wood2012power}.
Variations,
including minimum generation constraints
and the possibility of turning off generators,
are typically called the \emph{unit commitment problem}.
(See, \eg, \cite[Ch.~4]{wood2012power}
or the survey \cite{padhy2004unit}.)
These additions to the problem formulation,
which model the important limitations of many
types of power generation,
generally result in a nonconvex optimization problem.

\paragraph{\emph{Static optimal power flow.}}
The static optimal power flow problem
extends the optimal dispatch problem
by considering the spatial distribution of generators and loads
in a network.
In addition to planning operation of the generators,
the system operator must also consider how power flows
through this network to the loads.
It was first formulated in \cite{carpentier1962};
a modern treatment can be found in \cite[Ch.~8]{wood2012power}.
Good historical treatments of the development of optimization
for power systems can be found in \cite{happ1977optimal}
and \cite{cain2012history}.

\paragraph{\emph{DC optimal power flow.}}
Most formulations of optimal power flow consider AC power,
which typically results in a nonconvex problem.
This substantially complicates the formulation,
and we do not consider it in this paper.
Our formulation is similar to the so-called
\emph{DC optimal power flow},
or the \emph{network optimal power flow} problem described in
\cite{taylor2015convex}.
This simplified problem does not consider the physical
method by which power flows through the network,
and has the benefit of retaining convexity,
which we exploit in our exposition.
We also note that convexity raises
the possibility of a distributed solution method;
this idea is explored in \cite{kraning2014dynamic}.
The possibility of using a blockchain to coordinate transactions
in such a distributed method is considered in \cite{blockchain2017}, and a similar decentralized market
 structure is studied in \cite{liu2018novel}.

\paragraph{\emph{AC optimal power flow.}}
We do note that although the AC optimal power flow problem is not convex,
substantial progress has been made in the last 10 years
toward approximating the AC problem
using convex optimization.
These involve relaxing the (nonconvex) quadratic equality constraints
associated with AC power flow,
and result in second-order cone programs or semidefinite programs;
see \cite{lavaei2012zero} and \cite[Ch.~3]{taylor2015convex} for details.

\paragraph{\emph{Symbolic languages for optimization.}}
Object-oriented programming is ideal for software that encapsulates
technical details and offers a simple interface to (even advanced)
users. This has been used to develop languages for specifying 
optimization problems \cite{yalmip, cvx, cvxpy, convexjl, cvxr}. 
On top of these, 
domain-specific languages have been developed, for example 
for portfolio management in finance \cite{BBDKKNS:17}.

\subsection{Outline}
We start with a simple network power flow model,
and increase the complexity of the formulation in each subsequent section,
adding additional levels of complexity.
In \S\ref{c-stat-opf},
we begin with a basic network model,
representing the distribution of devices across a network.
This allows our formulation to capture spatial phenomena,
and in particular, the fact that the price of power can vary 
at different locations on a network.
(For example, power is typically cheaper close to cheap generators,
and more expensive close to loads,
especially if transmission is difficult or constrained.)
In \S\ref{c-dyn-opf},
we extend this model to account for phenomena occurring over time,
such as time-varying loads and availability of renewable power 
generation,
energy storage, generator ramp rate restrictions,
and deferrable loads.
Here we see that the price of power varies both across the network
and in time.
In \S\ref{c-mpc}, we use the dynamic formulation
for model predictive control, a method that replaces uncertain
future quantities with forecasts.
As seen in many other applications of model predictive control,
the feedback inherent in such a system gives good performance
even when the forecasts are not particularly accurate.
Finally, in \S\ref{c-robust-mpc},
we add an explicit uncertainty model to account for our prediction
or forecast errors, leading to an improved model predictive control
formulation.
We also present, in Appendix~\ref{c-appendix-forecasts}, a simple method for forecasting dynamic quantities, such as power availability
of renewable generators.

\section{Static optimal power flow}
\label{c-stat-opf}

In this section we describe our basic abstractions, which we
use throughout the paper (sometimes in more sophisticated forms).
Our abstractions follow \cite{kraning2014dynamic}.

In this section we work in a \emph{static} setting, \ie, we consider power flows that 
are constant over time, or more realistically, constant over some given
time interval such as one minute, 15 minutes, or one hour.
Thus any power that we refer to in this section can be converted to energy
by multiplying by the given time interval.

\subsection{Network model}
\label{s-network-model}

We work with three abstractions: \emph{devices}, \emph{terminals}, and \emph{nets}.
We first describe the setup in words, without equations;
after that, we introduce our formal notation.

\subsubsection{Devices and terminals}
\emph{Devices} produce, consume, and transfer power. 
Examples include generators, loads, transmission lines, and power converters.
Each device has one or more \emph{terminals},
across which power can flow (in either direction).
We adopt the sign convention that positive terminal power means power flows
into the device at that terminal; negative power corresponds to power flowing out
of the device at that terminal.
For example, we would expect a load (with a single terminal) to have positive power,
whereas a generator would have negative power at its (single) terminal.
As another example, a transmission line (or other energy transport or conversion device)
has two terminals; we would typically expect one of its terminal powers to be positive
(\ie, the terminal at which the power enters) and the other terminal power 
(at which the power leaves) to be negative.
(The sum of these two terminal powers is the net power entering the device,
which can be interpreted as the power lost or dissipated.)

We do not specify the physical transport mechanism by which power flows across terminals;
it is simply a power, measured in Watts (or kW, MW, or GW).
The physical transport could be a DC connection (at some specific
voltage), or a single- or multi-phase AC connection (at some specific voltage).
We do not model AC quantities like voltage magnitude, phase angle,
or reactive power flow.
In addition, power can have a different physical transport mechanism at its different terminals.
For example, the two terminals (in our sense, not the electrical sense) 
of an AC transformer transfer power at different voltages;
but we keep track only of the (real) power flow on the primary and secondary terminals.

Each device has a \emph{cost function},
which associates a (scalar) cost with its terminal powers.
This cost function can be used to model operating cost (say, of a generator), or 
amortized acquisition or maintenance cost.  
The cost can be infinite for some device terminal 
powers; we interpret this as indicating that the terminal powers violate a constraint or are 
impossible or infeasible for the device.
The cost function is a quantity that we would like to be small.  The negative of the 
cost function, called the \emph{utility function}, is a quantity that we would like to
be large.

\subsubsection{Nets}

\emph{Nets} exchange power between terminals.  
A net consists of a set of two or more
terminals (each of which is attached to a device).  If a terminal is in a net, we say
it is \emph{connected}, \emph{attached}, or \emph{adjacent} to the net.
At each net we have (perfect) power flow conservation; in other words,
the sum of the attached terminal powers is zero.  This means that
the sum of total power flowing from the net to the device terminals 
exactly balances the total power flowing into the net from device terminals.
A net imposes no constraints on the attached terminal powers other than conservation,
\ie, they sum to zero.
We can think of a net as an ideal bus with no power loss or limits imposed,
and without electrical details such as voltage, current, or AC phase angle.

A one-terminal net is not very interesting, since power conservation requires that the 
single attached terminal power is zero.
The smallest interesting net is a two-terminal net.  
The powers of the two connected terminals sum to zero; \ie, 
one is the negative of the other.  We can think of a two-terminal 
net as an ideal lossless power transfer point between two terminals;
the power flows from one terminal to the other.

\subsubsection{Network}
A \emph{network} is formed from a collection of devices and nets by connecting each terminal 
of each device to one of the nets.
The \emph{total cost} associated with the network is the sum of the costs of its devices,
which is a function of the device terminal powers.
We say that the set of terminal powers in a network is \emph{feasible} if the cost is finite,
and if power conservation holds at each net.
We say that the set of terminal powers in a network is \emph{optimal} if it is feasible,
and it minimizes the total cost among all feasible terminal power flows.
This concept of optimal terminal powers is the central one in this paper.

\subsubsection{Notation}
We now describe our notation for the abstractions introduced above.
We will use this notation (with some extensions described later) throughout this paper.

There are $D$ devices, indexed as $d=1,\ldots,D$.
Device $d$ has $M_d$ terminals,
and there are $M$ terminals in total
(\ie, $\sum_{d=1}^D M_d = M$).
We index terminals using $m=1, \ldots, M$.
We refer to this ordering of terminals as the \emph{global ordering}.
The set of all terminal powers is represented as a vector $p\in \reals^M$,
with $p_m$ the power flow on terminal $m$ (in the global ordering).
We refer to $p$ as the \emph{global power vector}; 
it describes all the power flows in the network.

The $M_d$ terminal powers of a specific device $d$ are 
denoted $p_d\in\reals^{M_d}$.
This involves a slight abuse of notation; we use $p_m$ to denote the (scalar)
power flow on terminal $m$ (under the global ordering); we use $p_d$ to 
denote the vector of terminal powers for the terminals of device $d$.
We refer to the scalar $(p_d)_i$ as the power on terminal $i$ of device $d$.
We refer to the ordering of terminals on a multi-terminal device as the
\emph{local ordering}.
For a single-terminal device, $p_d$ is a number (\ie, in $\reals$).

Each device power vector $p_d$ consists of a subvector or selection from the 
entries of the global power vector $p$.
We can express this as $p_d = B_d p$, where $B_d$ is the matrix that maps the
global terminal ordering into the terminal ordering for device $d$.
These matrices have the simple form
\[
  (B_d)_{ij}=
  \begin{cases}
    1 & \text{(global) terminal $j$ is the $i$th (local) terminal of device $d$}\\
    0 & \text{otherwise}.
  \end{cases}
\]
We refer to $B_d$ as the \emph{global-local matrix}, since it maps the global power vector
into the local device terminal powers.
For a single-terminal device, $B_d$ is a row vector, $e_k^T$, where $e_k$
is the $k$th standard unit vector, and $k$ is the global ordering index of 
the terminal.
The cost function for device $d$ is given by
$f_d:\reals^{M_d}\to \reals \cup \{\infty\}$.
The cost for device $d$ is 
\[
f_d(p_d)=f_d(B_dp).
\]

The $N$ nets are labeled $n=1,\ldots,N$.
Net $n$ contains $M_n$ terminals, and
we denote by $p_n$ the vector of powers
corresponding to the terminals in $n$, ordered under the global ordering.
(Because each terminal appears in exactly one net,
we have $\sum_{n=1}^N M_n = M$.)
Here too we abuse notation: $M_d$ is the number of terminals of device $d$,
whereas $M_n$ is the number of terminals in net $n$.
The symbol $p$ by itself always refers to the global power vector.
It can have two meanings when subscripted: $p_m$ is the (scalar)
power flow on (global)
terminal $m$; $p_d$ is the vector of power flows for device $d$.
The power flow on (local) terminal $i$ on device $d$ is
$(p_d)_i$.

The terminals in each net can be described by an \emph{adjacency matrix} 
$A \in \reals^{N \times M}$, defined as
\[
  A_{nm} =
  \begin{cases}
    1 &  \text{terminal $m$ is connected to net $n$}\\
    0 & \text{otherwise}.
  \end{cases}
\]
Each column of $A$ is a unit vector corresponding to a terminal; each 
row of $A$ corresponds to a net, and consists of a row vector with entries $1$ and $0$ 
indicating which nets are adjacent to it.
We will assume that every net has at least one adjacent terminal, so
every unit vector appears among the columns of $A$, which implies it is full rank.

The number $(Ap)_n$ is the sum of the terminal powers over terminals in net $n$,
so the $n$-vector $Ap$ gives the total or net power flow out of each net.
Conservation of power at the nets can then be expressed as
\begin{equation}
\label{e-power-conservation}
  Ap = 0,
\end{equation}
which are $n$ equalities.

The total cost of the network, denoted $f:\reals^M\to \reals$,
maps the power vector $p$ to the (scalar) cost.
It is the sum of all device costs in the network:
\begin{equation}
	f(p) = \sum_{d=1}^D f_d(p_d) = \sum_{d=1}^D f_d(B_d p).
  \label{e-objective-def}
\end{equation}
A power flow vector $p$ is called \emph{feasible} if $Ap=0$ and $f(p)< \infty$.
It is called \emph{optimal} if it is feasible, and has smallest cost among 
all feasible power flows.

\subsubsection{Example}
\begin{figure}
\begin{center}
\scalebox{.8}{\begin{tikzpicture}

\def\termsep{2cm}

\node (n1)  [net] {\,net 1\,};

\node (t1)   [device, position=240:{\termsep} from n1] {line 1};
\node (t2)   [device, position=300:{\termsep} from n1] {line 2};

\node (n2)  [net, position=240:{\termsep} from t1] {\,net 2\,};
\node (n3)  [net, position=300:{\termsep} from t2] {\,net 3\,};

\node (t3)   [device, distance={\termsep}, right of=n2] {line 3};
\node (g1)   [device, left of=n1] {gen.\ 1};
\node (g2)   [device, below of=n3] {gen.\ 2};
\node (l1) [device, right of=n1] {load 1};
\node (l2) [device, below of=n2] {load 2};

\path[terminal] (n1) edge[pos=.50, above] node {\textcolor{flow}{$p_{1 }$}}  (g1);
\path[terminal] (n1) edge[pos=.50, above] node {\textcolor{flow}{$p_{2 }$}}  (l1);
\path[terminal] (n1) edge[pos=.50, left ] node {\textcolor{flow}{$p_{3 }$}}  (t1);
\path[terminal] (n1) edge[pos=.50, left ] node {\textcolor{flow}{$p_{4 }$}}  (t2);
\path[terminal] (n2) edge[pos=.50, left ] node {\textcolor{flow}{$p_{5 }$}}  (t1);
\path[terminal] (n3) edge[pos=.50, left ] node {\textcolor{flow}{$p_{6 }$}}  (t2);
\path[terminal] (n2) edge[pos=.50, below] node {\textcolor{flow}{$p_{7 }$}}  (t3);
\path[terminal] (n3) edge[pos=.50, below] node {\textcolor{flow}{$p_{8 }$}}  (t3);
\path[terminal] (n2) edge[pos=.50, right] node {\textcolor{flow}{$p_{9 }$}}  (l2);
\path[terminal] (n3) edge[pos=.50, left ] node {\textcolor{flow}{$p_{10}$}}  (g2);
\end{tikzpicture}}
\qquad\qquad
\scalebox{.8}{\begin{tikzpicture}

\node (g1)  [net] {G1};

\draw [termnoarrow] (g1) 
-- ++(0cm, -1cm) node (b1r) {}
-- ++(0cm, -.5cm)
-- ++(5cm, 0cm)
-- ++(0cm, .5cm) node (b2l) {}
++(0cm, 1cm) node (g2) [net] {G2};

\draw [terminal] (g2) -- (b2l.center);

\draw (b1r) ++(-1cm, 0cm) node (b1l) {};
\draw [bus] (b1l) -- ++(-.5cm, 0cm) -- (b1r.center) -- ++(.5cm, 0cm);

\draw (b2l) ++(1cm, 0cm) node (b2r) {};
\draw [bus] (b2l.center) -- ++(-.5cm, 0cm) -- (b2r.center) -- ++(.5cm, 0cm);

\draw [terminal] (b1l.center) -- ++(0cm, 1cm);

\draw [terminal] (b1l.center) 
-- ++(0cm, -.5cm)
-- ++(3cm, -2.5cm)
-- ++(0cm, -.5cm) node (b3l) {};

\draw [termnoarrow] (b3l) ++(1cm, 0cm) node (b3r) {};
\draw [bus] (b3l.center) -- ++(-.5cm, 0cm) -- (b3r.center) -- ++(.5cm, 0cm);

\draw [terminal] (b3l) ++(.5cm, 0cm) -- ++(0cm, -1cm);

\draw [terminal] (b2r.center) 
-- ++(0cm, -.5cm)
-- ++(-3cm, -2.5cm)
-- ++(0cm, -.5cm) node (b3l) {};

\draw (b3l) ++(0cm, -2cm) node (bottom) {};

\end{tikzpicture}}
\caption{
\emph{Left.} The three-bus example network.
\emph{Right.} Traditional schematic.
    }
\label{f-transmission-network}
\end{center}
\end{figure}
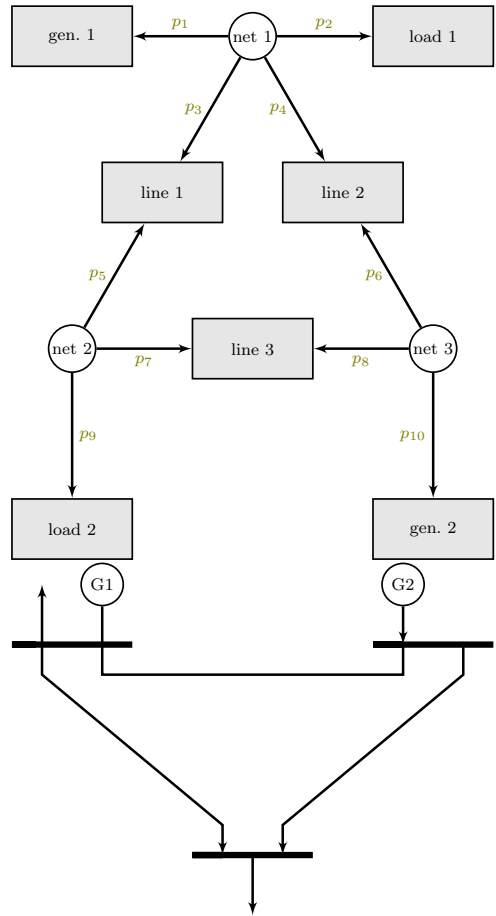

As an example of our framework,
consider the three-bus network shown
in figure~\ref{f-transmission-network}.
The two generators and two loads
are each represented as single-terminal devices,
while the three transmission lines,
which connect the three buses,
are each represented as two-terminal devices, so this
network has $D=7$ devices and $M=10$ terminals.
The three nets,
which are the connection points of these seven devices,
are represented in the figure as circles.
The device terminals are represented as lines (\ie, edges)
connecting a device and a net.
Note that our framework puts transmission lines
(and other power-transfer devices)
on an equal footing with other devices such as
generators and loads.

In the figure we have labeled the terminal powers with the global index.
For the network, we have
\[
A = 
\begin{bmatrix}
  1 & 1 & 1 & 1 & 0 & 0 & 0 & 0 & 0 & 0 \\
  0 & 0 & 0 & 0 & 1 & 0 & 1 & 0 & 1 & 0 \\
  0 & 0 & 0 & 0 & 0 & 1 & 0 & 1 & 0 & 1 
\end{bmatrix}.
\]
The conservation of power condition $Ap=0$ can be written explicitly as
\[
p_1 + p_2 + p_3 + p_4 = 0,
\qquad
p_5 + p_7 + p_9 = 0,
\qquad
p_6 + p_8 + p_{10} = 0
\]
(for nets~1, 2, and 3, respectively).
The third device is line 1, with global-local matrix
\[
B_3 = 
\begin{bmatrix}
  0 & 0 & 1 & 0 & 0 & 0 & 0 & 0 & 0 & 0 \\
  0 & 0 & 0 & 0 & 1 & 0 & 0 & 0 & 0 & 0 \\
\end{bmatrix}.
\]
The generators typically produce power,
not consume it, so we expect the generator powers $p_1$ and $p_{10}$ to be negative.
Similarly, we expect the load powers $p_2$ and $p_9$ to be positive,
since loads typically consume power.
If line 3 is lossless, we have $p_7 + p_8 = 0$;
if power is lost or dissipated in line 3, $p_7 + p_8$
(which is power lost or dissipated) is positive.

\subsubsection{Generators, loads, and transmission lines}
Our framework can model a very general network, with devices that have more than two
terminals, and devices that can either generate or consume power.  But here we describe 
a common situation, in which the devices fall into three general categories:  
Loads are single-terminal devices that consume power, \ie, have positive terminal
power. 
Generators are single-terminal devices that generate power, \ie, have negative
terminal power.  And finally, transmission lines and power conversion devices are 
two-terminal devices that transport power, possibly with dissipation, \ie, 
the sum of their two terminal powers is nonnegative.

For such a network, power conservation allows us to make a statement 
about aggregate powers.
Each net has total power zero, so summing over all nets we conclude that the sum of all
terminal powers is zero.  (This statement holds for any network.)
Now we partition the terminals into those associated with generators, those associated
with loads, and those associated with transmission lines.  Summing the terminal
powers over these three groups, we obtain the total generator power, the total load
power, and the total power dissipated or lost in the transmission lines.
These three powers add to zero.
The total generator power is negative, the total load power is positive, and
the total power dissipated in transmission lines is nonnegative.
Thus we find that the total power generated (expressed as a positive number)
exactly balances the total load, plus the total power loss in transmission lines.

\subsection{Optimal power flow}
\label{s-static-OPF}
The \emph{static optimal power flow problem}
consists of finding the terminal powers
that minimize the total cost of the network over
all feasible terminal powers:
\begin{equation}
  \begin{array}{ll}
    \mbox{minimize} & f(p) \\
    \mbox{subject to} & Ap = 0.
  \end{array}
  \label{e-static-prob}
\end{equation}
The decision variable is $p \in \reals^M$, the vector of all terminal powers.
The problem is specified by the cost
functions $f_d$ of the $D$ devices,
the adjacency matrix $A$,
and the global-local matrices $B_d$, for $d=1,\ldots,D$.
We refer to this problem as the static OPF problem.
We will let $p^\star$ denote an optimal power flow vector, and 
we refer to $f(p^\star)$ as the optimal cost for the power flow
problem~(\ref{e-static-prob}).
The OPF problem is a convex optimization problem 
if all the device cost functions are convex \cite{cvxbook}.
Roughly speaking, this means that it can be solved exactly and efficiently, even
for large networks.

\subsubsection{Optimality conditions}
\label{s-static-opt-conds}
If all the device cost functions are convex and differentiable,
a terminal power vector $p^\star \in \reals^M$ is optimal
for \eqref{e-static-prob}
if and only if there exists
a Lagrange multiplier vector $\lambda \in \reals^N$
such that
\begin{equation}
\nabla f(p^\star) = A^T\lambda, \qquad
Ap^\star = 0,
\label{e-static-opt-cond}
\end{equation}
where $\nabla f(p^\star)$ is the gradient of $f$ at $p^\star$
\cite{bertsekas2016nonlinear}.
The second equation is the conservation of power constraint
of the OPF problem (\ref{e-static-prob}).
For a given optimal flow vector $p^\star$,
there is a unique Lagrange multiplier vector $\lambda$
satisfying (\ref{e-static-opt-cond}).
(This follows since the matrix $A$ has full rank.)
The Lagrange multiplier vector $\lambda$ will come up again in
\S\ref{s-static-LMP}, where it will be interpreted as a vector of prices.

Some of the aforementioned assumptions
(convexity and differentiability of the cost function)
can be relaxed.
If the cost function is convex but not differentiable,
the optimality conditions~(\ref{e-static-opt-cond})
can be extended in a straightforward manner
by replacing the gradient with a subgradient.
(In this case, the Lagrange multiplier vector may not be unique.)
For a detailed discussion, see
\cite[\S 28]{rockafellar1997convex}.
If the cost function is differentiable but not convex,
the conditions~(\ref{e-static-opt-cond}) are necessary,
but not sufficient, for optimality;
see \cite[Ch.~4]{bertsekas2016nonlinear}.
When the cost function is neither convex nor differentiable, optimality 
conditions similar to~(\ref{e-static-opt-cond}) can be formulated using
generalized (Clarke) derivatives \cite{clarke1975generalized}.

\subsubsection{Solving the optimal power flow problem}
\label{s-static-convexity}
When all the device cost functions are convex, the objective function $f$ is
convex, and the OPF problem is a convex optimization problem. 
It can be solved exactly (and efficiently) using standard algorithms; see \cite{cvxbook};
all such methods also compute the Lagrange multiplier $\lambda$ as well
as an optimal power flow $p^\star$.

If any of the device cost functions is not convex, the OPF
problem is a nonconvex optimization problem.
In practical terms, this means that finding
(and certifying) a global solution to \eqref{e-static-prob}
is difficult in general.
Local optimization methods, however, can efficiently find 
power flows and a Lagrange multiplier vector that satisfy
the optimality conditions~(\ref{e-static-opt-cond}).


\subsection{Prices and payments}
\label{s-static-LMP}
In this section we describe a fundamental concept in power flow 
optimization, \emph{locational marginal prices}.  These prices lead
to a natural scheme for payments among the devices.

\subsubsection{Perturbed problem}
Suppose a network has optimal power flow $p^\star$,
and we imagine extracting additional power from each net.
We denote this perturbation by a vector $\delta\in\reals^N$. 
When $\delta_n>0$, we extract additional power from net $n$;
$\delta_n<0$ means we inject additional power into net $n$.
Taking these additional power flows into account,
the power conservation constraint $Ap=0$ becomes $Ap + \delta = 0$.
The \emph{perturbed optimal power flow problem} is then 
\[
	\begin{array}{ll}
    \mbox{minimize} & f(p) \\
    \mbox{subject to} & Ap  + \delta = 0.
  \end{array}
\]
Note that when $\delta=0$, this reduces to the optimal power flow problem.

We define $F:\reals^N \to \reals \cup \{\infty\}$, the \emph{perturbed
optimal cost function}, as the optimal cost of the perturbed optimal
power flow problem, which is a function of $\delta$.
Roughly speaking, $F(\delta)$ is the minimum total network cost, obtained by 
optimizing over all network power flows, taking into account the net power
perturbation $\delta$.
We can have $F(\delta)=\infty$, which means that with the perturbed
power injections and extractions,
there is no feasible power flow for the network.
The optimal cost of the unperturbed network is $F(0)$.

\subsubsection{Prices}
The change in the optimal cost from the unperturbed 
network is $F(\delta)-F(0)$.
Now suppose that $F$ is differentiable at $0$ 
(which it need not be; we discuss this below).
We can approximate the cost change, for small perturbations, as
\[
F(\delta) - F(0) \approx \nabla F(0)^T \delta = \sum_{n=1}^N
\frac{\partial F}{\partial \delta_n}(0) \delta_n.
\]
This shows that the approximate change in optimal cost 
is a sum of terms, each 
associated with one net and proportional to the perturbation power $\delta_n$.
We define the \emph{locational marginal price} (or just \emph{price}) at net $n$ to be
\[
\frac{\partial F}{\partial \delta_n}(0).
\]

The locational marginal price at net $n$ has a simple interpretation.  
We imagine a network operating at an optimal power vector $p^\star$.
Then we imagine that a small amount of additional power is extracted 
from net $n$.  We now re-optimize all the power flows, taking into account
this additional power perturbation.  The new optimal cost will (typically) 
rise a bit from the unperturbed value, by an amount very close to
the size of our perturbation times the locational marginal
price at net $n$.

It is a basic (and easily shown) result in optimization
that, when $f$ is convex and differentiable, and $F$ is differentiable,
we have
\cite{schweppe1988spot, papavasiliou2017analysis}
\begin{equation}
\lambda = \nabla F(0).
\label{e-price-def}
\end{equation}
In other words, the Lagrange multiplier in the OPF optimality 
condition~(\ref{e-static-opt-cond}) is precisely the vector of 
locational marginal prices.

Under usual circumstances, the prices are positive, which means that when we 
extract additional power from a net, the optimal system cost increases;
for example, at least one generator or other power provider must increase 
its power output to supply the additional power extracted, so its cost (typically)
increases.
In some pathological situations, locational marginal prices can be negative.
This means that by extracting power from the net, we can \emph{decrease}
the total system cost.  While this can happen in practice, we consider it to be 
a sign of poor network design or operation.

If $F$ is not differentiable at $0$, it is still possible
to define the prices, but the treatment becomes complicated and
mathematically intricate, so we do not include it here.
When the OPF problem is convex, the system cost function $F$ is convex, and
the prices would be given by a subgradient of $F$ at $0$;
see \cite[\S 28]{rockafellar1997convex}.
In this case, the prices need not be unique.
When the OPF problem is differentiable but not convex, the vector $\lambda$ in the 
(local) optimality condition (\ref{e-static-opt-cond}) can be interpreted
as predicting the change in local optimal cost with net perturbations.

\subsubsection{Payments}
\label{s-static-payments}
The locational marginal prices provide the basis for a natural payment scheme
among the devices.
With each terminal $i$ (in the global ordering) we associate a payment
(by its associated device) equal to its power times the associated net price.
We sum these payments over the terminals in each device to obtain the payment 
$P_d$ that is to be made by device $d$.
Define $\lambda_d$ as the vector of prices at the nets
containing the terminals of device $d$, \ie, 
\[
\lambda_d = B_dA^T \lambda.
\]
(These prices are given in the local terminal ordering for device $d$.)
The payment from device $d$
is the power flow at its terminals multiplied by the 
corresponding locational marginal prices, \ie,
\[
P_d = \lambda_d^Tp_d^\star.
\]

For a single-terminal device, this reduces to paying for the power consumed
(\ie, $p_d^\star$) at a rate given by the locational marginal price (\ie, $\lambda_d$).
A generator would typically have $p_d^\star<0$, and as mentioned above,
we typically have $\lambda_d>0$, so the payment is negative, \ie, it is
income to the generator.

For a two-terminal device, the payment is the sum of the two payments associated 
with each of the two terminals.
For a transmission line or other power transport or conversion device,
we typically have one terminal power positive (where power enters the 
device) and one terminal power negative (where power is delivered).
When the adjacent prices are positive, such a device receives payment for 
power delivered, and pays for the power where it enters.  The payment by the device is 
typically smaller than the payment to the device, so it typically
derives an income, which can be considered its compensation for transporting the power.

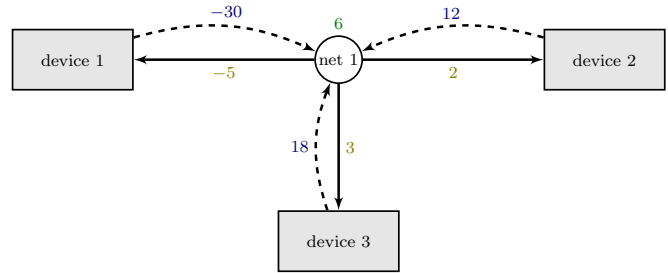
\begin{figure}
\begin{center}
\scalebox{.8}{\begin{tikzpicture}

\def\termsep{3cm}

\node (n1)  [net, label=above:{\textcolor{price}{$6$}}] {\,net 1\,};

\node (d1)   [device, position=180:{\termsep   } from n1] {device 1};
\node (d2)   [device, position=000:{\termsep   } from n1] {device 2};
\node (d3)   [device, position=270:{\termsep*.7} from n1] {device 3};

\path[terminal] (n1) edge[pos=.50, right] node {\textcolor{flow}{$3$}}  (d3);
\path[terminal] (n1) edge[pos=.50, below] node {\textcolor{flow}{$2$}}  (d2);
\path[terminal] (n1) edge[pos=.50, below] node {\textcolor{flow}{$-5$}}  (d1);

\path[terminal, dashed] (d3) edge[pos=.50, left, bend left=20] 
    node {\textcolor{pmnt}{$18$}}  (n1);
\path[terminal, dashed] (d2) edge[pos=.50, above, bend right=20]
    node {\textcolor{pmnt}{$12$}}  (n1);
\path[terminal, dashed] (d1) edge[pos=.50, above,  bend left=20] 
    node {\textcolor{pmnt}{$-30$}}  (n1);

\end{tikzpicture}}
\caption{
A net with 3 terminals. The dark arrows show the 
power flows (with values labeled in yellow),
and the dashed arrows show the three device payments (with values
labeled in blue).
The net price is shown in green.}
\label{f-net-payments}
\end{center}
\end{figure}

The total of all payments at net $n$ is $\lambda_n$ times the sum of the 
powers at the net.  But the latter is zero, by power conservation, so the total
of all payments by devices connected to a net is zero.
Thus the device payments can be thought of as an exchange of money taking place at 
the nets; just as power is preserved at nets, so are the payments.
We can think of nets as handling the transfer of power
among devices, as well as the transfer of money (\ie, payments), 
at the rate given by its locational marginal price.
This idea is illustrated in figure~\ref{f-net-payments}, where the dark lines
show power flow, and the dashed lines show payments, \ie, money flow.  Both 
are conserved at a net.
In this example, device~1 is a generator supplying power to devices~2 and~3, which
are loads. Each of the loads pays for their power at the locational marginal price;
the sum of the two payments is income to the generator.

Since the sum of all payments at each net is zero, it follows that the sum of all
payments by all devices in the network is zero.
(This is also seen directly: the sum of all device payments is $\lambda^T A p^\star$,
and we have $A p^\star = 0$.)
This means that the payment scheme is an exchange of money among the devices,
at the nets.

\subsubsection{Profit maximization}
\label{s-static-profits}
According to the payment scheme described above, 
device $d$ pays for power at its terminals
at rates given by the device price vector $\lambda_d$.
If we interpret the device cost function $f_d$ as a cost (in the same
units as the terminal payments), the device's net revenue or profit is
\[
  -\lambda_d^T p_d - f_d(p_d).
\]
We can think of the first term as the revenue associated with the power produced
or consumed at its terminals; the second term is the cost of operating the device.

When $f_d$ is differentiable, this profit is maximized when
\begin{equation}
\nabla f_d(p_d) + \lambda_d = 0.
\label{e-supply-fun}
\end{equation}
This is the first equation of \eqref{e-static-opt-cond}
(when evaluated at $p^\star$).
We conclude that, given the locational marginal prices $\lambda$,
the optimal power vector $p^\star$
maximizes each device's individual profit.
(This assumes that device $d$ acts as a \emph{price taker},
\ie, it maximizes its profit
while disregarding the indirect impact of its terminal power flows
on the marginal prices of power at neighboring nets.
Violations of this assumption can result in deviations from optimality;
see \cite[Ch.~8]{luenberger1995microeconomic}
and \cite[\S6.3]{taylor2015convex}.)
The same profit maximization principle can be established when
$f_d$ is convex but not differentiable.  In this case any optimal
device power $p_d$ maximizes the profit
$-\lambda_d^T p_d - f_d(p_d)$.  
(But in this case, this does not determine the device optimal power uniquely.)

Note that (\ref{e-supply-fun}) relates the price of power
at adjacent nets to the (optimal) power consumed by the device.
For a single-terminal device with differentiable $f_d$
and optimal power $p_d^\star$, we see that the adjacent net price
must be $\lambda_d = -f_d'(p_d^\star)$.
This is the \emph{demand function} for the device.
When $f_d'$ is invertible, we obtain 
\BEQ\label{e-scalar-demand}
p_d^\star  = \left( f_d' \right) ^{-1} (-\lambda_d),
\EEQ
which can be interpreted as a prescription for how much power to consume or produce
as a function of the adjacent price.

For a multi-terminal device with differentiable $f_d$, and optimal power
$p_d^\star$, the vector of 
prices $\lambda_d$ at the adjacent nets is $\lambda_d = - \nabla f_d(p_d^\star)$,
which is the (multi-terminal) demand function for the device.
When the device gradient function is invertible, its inverse 
maps the (negative) adjacent net prices 
into the power generated or supplied by the device.
\BEQ\label{e-vector-demand}
p_d^\star = \left( \nabla f_d \right)^{-1} (-\lambda_d).
\EEQ
(In the case of nondifferentiable, convex cost functions,
the subgradient mapping is used here,
and the demand function is set valued.)

The demand functions or their inverses~(\ref{e-scalar-demand}) 
and~(\ref{e-vector-demand}) 
can be used to derive a suitable cost function for a device.
For example if a single-terminal device $d$ connected to net $n$ 
uses (decreasing, invertible) demand function $p_d = D_d(\lambda_n)$, 
we have $\lambda_n = D_d^{-1} (p_d) = -f_d'(p_d)$,
and we can take as cost function
\[
f_d(p_d) =  -\int_{0}^{p_d}D_d^{-1} (u)   \; du,
\]
which is convex.

The discussion above for single-terminal devices 
uses the language appropriate when the device 
represents a load, \ie, has positive terminal power, and $f_d$ is typically
decreasing.  While the equations still hold, the language would change 
when the device represents a generator, \ie, has negative terminal
power and is typically decreasing.

\subsection{Device examples}
\label{s-static-devices}
Here we list several practical device examples.
All cost functions discussed here are convex, unless otherwise noted.
We also discuss device constraints;
the meaning is that power flows that do not satisfy the device constraints 
result in infinite cost for that device.

\subsubsection{Generators} \label{s-static-generators}
A generator is a single-terminal device
that produces power,
\ie, its terminal power $p_d$ satisfies $p_d \leq 0$.
We interpret $-p_d$ as the (nonnegative) power generated.   
The device cost $f_d(p_d)$ is
the cost of generating power $-p_d$, and $-f_d'(p_d)$ is the marginal cost of
power when operating at power $p_d$.

A generic generator cost function has the form
\begin{equation}
f_d(p_d) =
\begin{cases}
\phi_d(-p_d) & p_{\rm min} \leq -p_d \leq p_{\rm max}  \\
\infty & \mbox{otherwise,}
\end{cases}
\label{e-generic-gen}
\end{equation}
where $p_{\rm min}$ and $p_{\rm max}$ are the minimum and maximum
possible values of generator power, and $\phi_d(u)$ is the cost of 
generating power $u$, which is convex and typically increasing.
Convexity means that the marginal cost of generating power is nondecreasing 
as the power generated increases.  When the generation cost is increasing, it means
that the generator `prefers' to generate less power.

The profit maximization principle connects the net price $\lambda_d$
to the generator power $p_d$.
When $-p_d^\star$ lies between $p_{\rm min}$ and $p_{\rm max}$, we have 
\[
\lambda_d = - f_d'(p_d^\star) = \phi'_d(-p_d^\star),
\]
\ie, the net price is the (nonnegative) marginal cost for the generator.
When $p_d^\star = p_{\rm min}$, we must have $\lambda_d \leq \phi_d'(p_{\rm min})$.
When $p_d^\star = p_{\rm max}$, we must have $\lambda_d \geq \phi_d'(p_{\rm max})$.
Since $\phi_d$ is convex, $\phi_d'$ is nondecreasing, so
$\phi_d'(p_{\rm min})$ and
$\phi_d'(p_{\rm max})$ are the minimum and maximum marginal costs for the 
generator, respectively.
Roughly speaking, the generator operates at its minimum power when the price is 
below the minimum marginal cost, and it operates at its maximum power when
the price is above its maximum marginal cost; when the price is in between,
the generator operates at a point where its marginal cost matches the net price.

\paragraph{\emph{Conventional quadratic generator.}}
A simple model of generator uses the generation cost function 
\begin{equation}
\phi_d(p_d) = 
\alpha p_d^2 + \beta p_d +\gamma,
\label{e-quadratic-gen}
\end{equation}
where $\alpha$, $\beta$, and $\gamma$ are parameters.
For convexity, we require $\alpha\geq 0$.
(We also typically have $\beta \leq 0$.) 
The value of the constant cost term $\gamma$ has no effect on the optimal 
power of the generator.

\paragraph{\emph{Fixed-power generator.}}
A fixed-power generator produces $p_{\rm fix}$ units of power;
this is an instance of the generic generator 
with $p_{\rm min} = p_{\rm max} = p_{\rm fix}$.  (The function $\phi$ is only 
defined for $u=p_{\rm fix}$, and its value has no effect on the optimal power,
so we can take it to be zero.)
A fixed-power generator places no constraint on the adjacent net price.

\paragraph{\emph{Renewable generator.}}
A renewable generator can provide any amount of power between
$0$ and $p_{\rm avail}$,
and does so at no cost,
where $p_{\rm avail} \geq 0$ is
the power available for generation.
It too is an instance of the generic generator, with $p_{\rm min}=0$,
$p_{\rm max}=p_{\rm avail}$, and $\phi(u)=0$.

The profit maximization principle tells us that if the adjacent net price
is positive, we have $p_d = p_{\rm avail}$;
if the adjacent net price is negative, we have $p_d = 0$.
In other words, a renewable generator operates (under optimality) at
its full available power if the net price is positive, and shuts down if it is 
negative.  If the generator operates at a power in between $0$ and $p_{\rm avail}$,
the adjacent price is zero.

\subsubsection{Loads}

A load is a single-terminal device that consumes power, \ie, $p_d \geq 0$.
We interpret $f_d(p_d)$ as the operating cost for consuming power $p_d$.
We can interpret $-f_d(p_d)$ as the utility to the load of consuming power $p_d$.
This cost is typically decreasing, \ie, loads `prefer' to consume more power.
The marginal utility is $-f'_d(p_d)$.
Convexity of $f_d$, which is the same as concavity of the utility,
means that the marginal utility of a load is nonincreasing with
increasing power consumed.

A generic load cost function has the form
\begin{equation}
f_d(p_d) =
\begin{cases}
\phi_d(p_d) & p_{\rm min} \leq p_d \leq p_{\rm max}  \\
\infty & \mbox{otherwise,}
\end{cases}
\label{e-generic-load}
\end{equation}
where $p_{\rm min}$ and $p_{\rm max}$ are the minimum and maximum
possible values of load power, and $\phi_d(u)$ is the cost of 
consuming power $u$, which is convex and typically decreasing.

The profit maximization principle connects the net price $\lambda_d$
to the load power $p_d$.
When $p_d^\star$ lies between $p_{\rm min}$ and $p_{\rm max}$, we have 
\[
\lambda_d = - f_d'(p_d^\star) = - \phi'_d(p_d^\star),
\]
\ie, the net price is the (typically nonnegative) marginal utility for the load.
When $p_d^\star = p_{\rm min}$, we must have $\lambda_d \geq - \phi_d'(p_{\rm min})$.
When $p_d^\star = p_{\rm max}$, we must have $\lambda_d \leq - \phi_d'(p_{\rm max})$.
Since $\phi_d$ is convex, $\phi_d'$ is nondecreasing, so
$-\phi_d'(p_{\rm min})$ and
$-\phi_d'(p_{\rm max})$ are the minimum and maximum marginal utilities for the 
load, respectively.
Roughly speaking, the load operates at its minimum power when the price is 
above the maximum marginal utility, and it operates at its maximum power when
the price is below its minimum marginal cost.

\paragraph{\emph{Fixed load.}}
A fixed load consumes a fixed amount $p_{\rm fix}>0$ of power;
\ie, the device power flow $p_d$ satisfies $p_d = p_{\rm fix}$.
It is an instance of the generic load with $p_{\rm min}=p_{\rm max} = p_{\rm fix}$.
The value of $f_d(p_{\rm fix})$ does not affect the power, so we can take it to be 
zero.

\paragraph{\emph{Power dissipation device.}}
A power dissipation device has no operating cost,
and can consume (dissipate) any nonnegative power.
This is an instance of our generic load, with $p_{\rm min} =0$,
$p_{\rm max}=\infty$, and for $p_d \geq 0$, $f_d(p_d)=0$.

\paragraph{\emph{Curtailable load.}}
A curtailable load has a desired or target power consumption level $p_{\rm des}$,
and a minimum allowable power consumption $p_{\rm min}$.
If it consumes less power than its desired value,
a penalty is incurred on the shortfall,
with a price $\lambda_{\rm curt} > 0$.
The cost is
\[
	f_d(p_d) =
	\begin{cases}
		\lambda_{\rm curt}(p_{\rm des} - p_d) & p_{\rm min} \leq p_d \leq 
p_{\rm des}\\
		\infty & \mbox{otherwise.}
	\end{cases}
\]
A curtailable load is also an instance of our generic load.

If the adjacent net price is less than $\lambda_{\rm curt}$, we have 
$p_d^\star = p_{\rm des}$, \ie, there is no shortfall.
If the adjacent net price it is greater than $\lambda_{\rm curt}$,
the load consumes its minimum possible power $p_{\rm min}$.
If the adjacent net price is $\lambda_{\rm curt}$, the curtailable 
device can consume any power between $p_{\rm min}$ and $p_{\rm des}$.

\subsubsection{Grid ties}
\label{s-grid-tie}
A grid tie is a single-terminal device representing a connection
to an external power grid.
When $p_d \geq 0$, we interpret it as power being injected into the grid.
When $p_d < 0$, we interpret $-p_d$ as the power extracted from the grid.

It is possible to buy power from the grid at price
$\lambda_{\rm buy}$,
and sell power to the grid at price
$\lambda_{\rm sell}$.
We assume arbitrage-free nonnegative prices,
\ie, $\lambda_{\rm buy} \geq \lambda_{\rm sell} \geq 0$.
From the point of view of our system optimizer,
the cost of the grid tie is the cost of 
power bought from (or sold to) the grid, \ie,
\[
f_d(p_d) = \max\{ -\lambda_{\rm buy} p_d, -\lambda_{\rm sell} p_d \}.
\]
(Recall that $-p_d$ is power that we take from the grid connection, so 
we are buying power when $-p_d>0$, and selling power when $-p_d<0$, \ie,
$p_d>0$.)

Any net adjacent to a grid tie with $p_d<0$ (\ie, power flows from the grid)
has price $\lambda_{\rm buy}$; when $p_d>0$ (\ie, power flows into the grid)
it has price $\lambda_{\rm sell}$.  When $p_d=0$, the adjacent net price 
is not determined, but must be between the buy and sell prices.

As a variation on the basic grid tie device,
we can add lower and upper limits on $-p_d$,
representing the maximum possible power we can sell or buy,
\[
-p_{\rm max,sell} \leq -p_d \leq -p_{\rm max,buy},
\]
where
$p_{\rm max,sell}\geq 0$ and $p_{\rm max,buy}\geq 0$ are the maximum power we 
can sell and buy, respectively.


\subsubsection{Transmission lines and converters}
\label{s-transmission-lines}
Here we consider
transmission lines, power converters, transformers,
and other two-terminal devices that transfer power
between their terminals,
As usual, constraints on the power flows
are encoded into the cost function.
We denote the power flow on 
distinct terminals of a device with numbered subscripts,
for example as $p_1$ and $p_2$.
(That is, the subscripts~1 and~2 are in the local device terminal ordering.)
We interpret $p_1+p_2$ as the net power consumed or dissipated 
by the device, which is typically nonnegative.
If $p_1>0$ and $p_2$ is negative, then power enters the device at terminal~1
and exits the device at terminal~2, and vice versa when $p_1<0$ and $p_2>0$.

\paragraph{\emph{Transmission lines.}}
An ideal, lossless transmission line (or power converter)
has a cost of zero,
provided that the power conservation constraint
\[
p_1 + p_2 = 0
\]
is satisfied,
where $p_1$ and $p_2$ are the power flows
into the two terminals of the device.
(Note that such an ideal power converter is the same as a two-terminal net.)

Additionally, we can enforce power limits
\[
p_{\rm min} \leq p_1 \leq p_{\rm max}.
\]
(This is the same as requiring
$-p_{\rm max} \leq p_2 \leq -p_{\rm min}$.)
The resulting cost function
(with or without power limits)
is convex.
When $p_{\rm min}=-p_{\rm max}$, the transmission line 
or converter is \emph{symmetric}, \ie, its cost function is 
the same if we swap $p_1$ and $p_2$.
When this is not the case, the device is \emph{directional} or 
\emph{oriented}; roughly
speaking, its terminals cannot be swapped.

For a lossless transmission line for which the limits are not active,
\ie, $p_{\rm min} < p_1 < p_{\rm max}$,
the prices at the two adjacent nets must be the same.
When the prices at the two adjacent nets are not the same,
the transmission line operates at its limit, with power flowing into the 
device at the lower priced net, and flowing out to the higher priced
net.   In this case the device is paid for transporting power.  For example
with $\lambda_1  < \lambda_2$, we have $p_1^\star = p_{\rm max}$ and 
$p_2^\star = -p_{\rm max}$, and the device earns 
revenue $p_{\rm max}(\lambda_2-\lambda_1)$.

\paragraph{\emph{Lossless transmission line with quadratic cost.}}
We can add a cost to a lossless transmission line, $\alpha p_1^2$, with
$\alpha >0$.  This (convex) cost discourages large power flows across 
the device, when the optimal power flow problem is solved.
This objective term alone does not model quadratic
losses in the transmission line, which would result in $p_1+p_2>0$; while 
it discourages power flow in the transmission line, it does not take into account
the power lost in transmission.
For such a transmission line, the 
profit maximization principle implies that power flows from the terminal
with higher price to the terminal of lower price, with a flow proportional
to the difference in price between the two terminals.

\begin{figure}
\begin{center}
\includegraphics[width=.9\columnwidth]{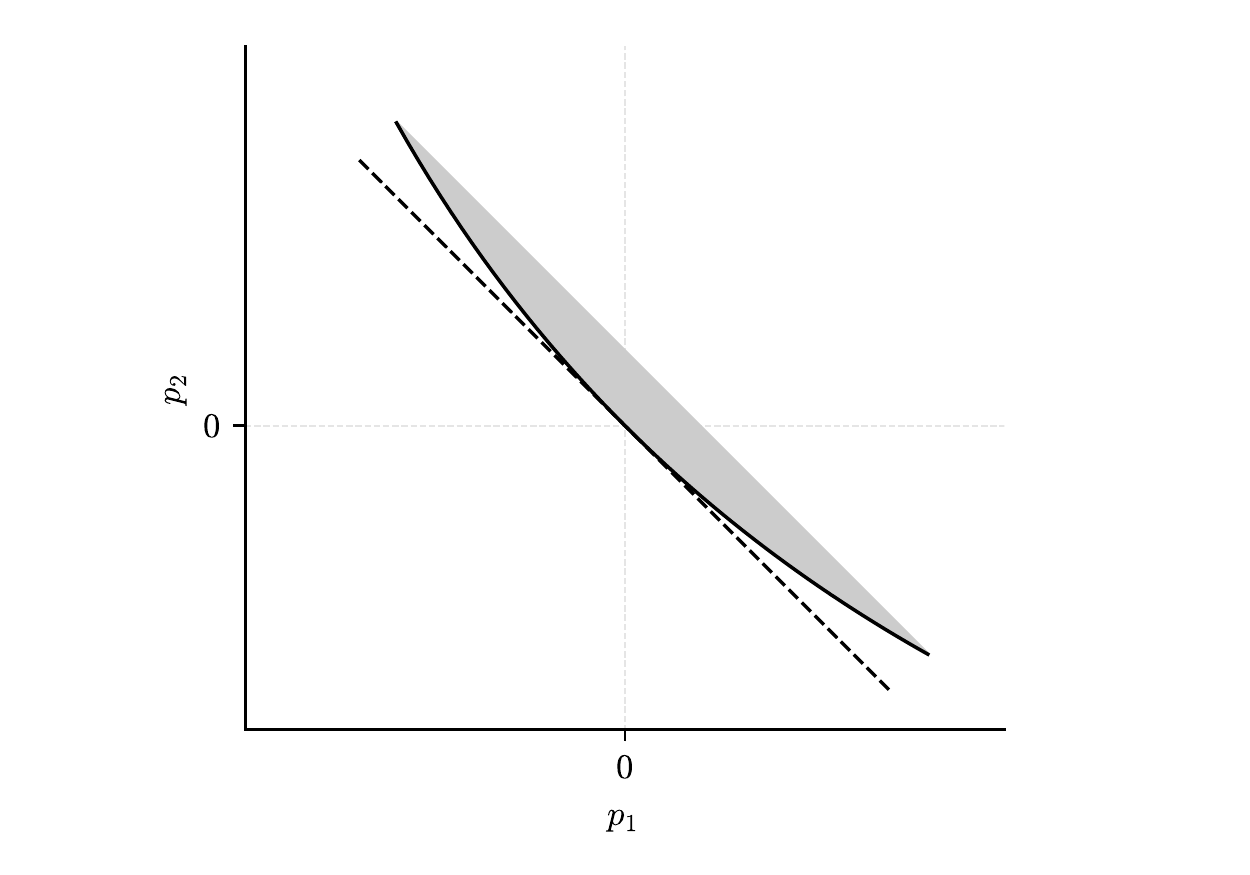}
\caption{
The set of feasible power flows, and its convex relaxation,
for a lossy transmission line.
}
\label{f-trans-line-cost}
\end{center}
\end{figure}

\paragraph{\emph{Transmission line with quadratic loss.}}
\label{par-lossy-transmission-line}
We consider a bi-directional transmission line 
with power loss $\alpha ((p_1 - p_2)/2)^2$
with parameter $\alpha>0$, and limit $|(p_1-p_2)/2| \leq p_{\rm max}$.
(This loss model can be interpreted as an
approximation for resistive loss.)
This model has constraints
\begin{align}
\label{e-lossy-trans-line}
p_1 + p_2 = \alpha \left( \frac{p_1 - p_2}{2}\right)^2, 
\quad \left|\frac{p_1-p_2}{2}\right| \leq p_{\rm max}.
\end{align}
The set of powers $p_1$ and $p_2$
that satisfy the above conditions
is the dark line shown in figure~\ref{f-trans-line-cost}.
Note that this set (and thus the device cost function)
is \emph{not} convex.
The ends of the curve are at the points
\[
p_1 = \alpha p_{\rm max}^2/2 + p_{\rm max}, \quad
p_2 = \alpha p_{\rm max}^2/2- p_{\rm max},
\]
and
\[
p_1 = \alpha p_{\rm max}^2/2 - p_{\rm max}, \quad
p_2 = \alpha p_{\rm max}^2/2 + p_{\rm max}.
\]

One way to retain convexity is to approximate 
this set with its convex hull, described by the constraints
\begin{align*}
\alpha p_{\rm max}^2 \geq p_1 + p_2 \geq \alpha \left(\frac{p_1 - p_2}{2}\right)^2.
\end{align*}
With this approximation,
the set of feasible powers is
the shaded region in figure~\ref{f-trans-line-cost}.
Power flows that are in this region, but not on the dark line,
correspond to artificially wasting power.
This approximation is provably exact
(\ie, the power flows $(p_1,p_2)$ lie on the dark line)
if the optimal price at a neighboring net is positive.
In practice, we expect this condition to hold in most cases.

\paragraph{\emph{Constant-efficiency converter.}}
We consider terminal~1 as the input and terminal~2 as the output of a converter
(in forward mode).
The forward conversion efficiency is given by $\eta \in (0,1)$,
and the reverse conversion efficiency is given by $\tilde \eta \in (0,1)$.
The device is characterized by
\begin{align*}
p_2 = \max \{-\eta p_1, -(1/\tilde \eta) p_1 \}, \qquad p_{\rm min} \leq p_1 \leq
p_{\rm max},
\end{align*}
Here $p_{\rm max}\geq 0$ is the maximum input power to the converter; 
$-p_{\rm min} \geq 0$ is the maximum power that can be extracted from the
converter.
For $p_1 \geq 0$ (forward conversion mode), the conversion loss is 
$p_1+p_2 = (1-\eta) p_1$, justifying the interpretation
of $\eta$ as the forward efficiency.
For $p_1 \leq 0$ (reverse conversion mode, which implies $p_2 \geq 0$),
the loss is $p_1+p_2 = (1-\tilde \eta) p_2$.

This set of constraints (and therefore the cost function) is \emph{not} convex,
but we can form a convex relaxation
(similar to the case of a lossy transmission line) which in this case is described
by a triangular region.
The approximation is exact if prices are positive at an adjacent net.

\subsubsection{Composite devices}
\label{s-composite-device}
\begin{figure}
\begin{center}
\begin{tikzpicture}

\node (net)  [net] {};
\node (d1)  [device, position=90:{1cm} from net] {$f_1$};
\node (d2)  [device, position=270:{1cm} from net] {$f_2$};
\node (dX)  [position=000:{2cm} from net] {\textcolor{flow}{$p_d$}};
\node (label)  at ($(d1.north west)+(-0.0cm,0.6cm)$) {$f_d$};

\node (rightend) [position=000:{3cm} from net] {};
\node (leftend)  [position=180:{3cm} from net] {};

\path[terminal] (net) edge[pos=.50, left] node {\textcolor{flow}{$p_1$}}  (d1);
\path[terminal] (net) edge[pos=.50, left] node {\textcolor{flow}{$p_2$}}  (d2);
\path[terminal] node {} (dX) edge[] (net);

\draw[dashed] ($(d1.north west)+(-0.3cm,0.3cm)$)  rectangle
              ($(d2.south east)+(0.3cm,-0.3cm)$);

\end{tikzpicture}
\caption{
A composite device can be formed
from two or more other devices.
}
\label{f-composite-device}
\end{center}
\end{figure}
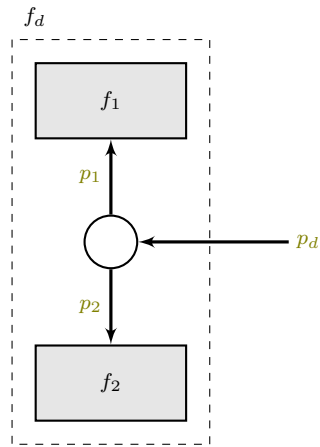

The interconnection of several devices at a single net
can itself be modeled as a single composite device.
We illustrate this with a composite device consisting of two 
single-terminal devices connected to a net with three terminals, one
of which is external and forms the single terminal of the 
composite device.
(See figure~\ref{f-composite-device}.)
For this example we have composite device cost function
\BEQ\label{e-min-conv}
f_d(p_d)
= \min _{\tilde p_1 + \tilde p_2 = p_d} f_1(\tilde p_1) + f_2(\tilde p_2).
\EEQ
(The function $f_d$ is the \emph{infimal convolution} of the functions
$f_1$ and $f_2$; see \cite[\S 5]{rockafellar1997convex}.)
This composite device can be connected to any network, and the optimal 
power flows will match the optimal power flows when the device is replaced
by the subnetwork consisting of the two devices and extra net.
The composite cost function~(\ref{e-min-conv}) is easily interpreted:  Given 
an external power flow $p_d$ into the composite device, it splits into the
powers $p_1$ and $p_2$ (as the net requires) in such a way as to minimize the 
sum of the costs of the two internal devices.
The composite device cost function~(\ref{e-min-conv}) is convex if the two 
component device cost functions are convex.

Composite devices can be also formed
from other, more complicated networks of devices, and can expose multiple external
terminals.
The composite device cost function in such cases is a simple generalization of
the infimal convolution~(\ref{e-min-conv}) for the simple case of two 
internal devices and one external terminal.
Such composite devices preserve convexity: Any composite device formed from
devices with convex cost functions also has a convex cost function.

Composite devices can simplify modeling.  For example, a wind farm, solar array,
local storage, and a transmission line that connects them to a larger grid can be 
modeled as one device.
We also note that there is no need to analytically compute the composite 
device function for use in a modeling system.  Instead we simply introduce the
sum of the internal cost functions, along with additional variables representing
the internal power flows, and the constraints associated with internal nets.
(This is the same technique used in the convex optimization modeling systems
\texttt{CVX} \cite{gb08,cvx} and \texttt{CVXPY} \cite{cvxpy} to represent compositions
of convex functions.)


\subsection{Network examples}
\label{sec-ex-load-gen}

\subsubsection{Two-device example}
\label{sec-ex-two-devices}
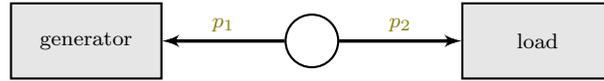
\begin{figure}
\begin{center}
  \scalebox{1.0}{\begin{tikzpicture}
\node (gen)  [device] {generator};
\node (net)  [net, right of=gen] {};
\node (load) [device, right of=net] {load};

\path [terminal] (net) -- (gen);
\path [terminal] (net) -- (load);

\path[terminal] (net) edge[pos=.50, above] node {\textcolor{flow}{$p_{1 }$}}  (gen);
\path[terminal] (net) edge[pos=.50, above] node {\textcolor{flow}{$p_{2 }$}}  (load);
\end{tikzpicture}}
  \caption{A network with one generator and one load connected to one net.}
  \label{f-gen_load_network}
\end{center}
\end{figure}

We consider the case of a generator
and a load connected to a single net,
as shown in figure~\ref{f-gen_load_network}.
For this network topology, the static OPF problem 
\eqref{e-static-prob} is
\begin{equation}
  \begin{array}{ll}
    \mbox{minimize} & f_{\rm gen}(p_1) + f_{\rm  load}(p_2) \\
    \mbox{subject to} & p_1 + p_2 = 0.
  \end{array}
  \label{e-example-prob}
\end{equation}
Assuming the cost function is differentiable at $p^\star$,
the optimality condition is
\begin{align*}
  f_{\rm gen}'(p_1^\star) = \lambda,\qquad
  f_{\rm load}'(p_2^\star) = \lambda,\qquad
  p_1^\star + p_2^\star = 0.
\end{align*}
We interpret $\lambda$ as the price, and
$f_{\rm gen}'(p_1^\star)$ and
$f_{\rm load}'(p_2^\star)$ as the marginal costs of the generator 
and load, respectively.

\begin{figure}
\begin{center}
  \includegraphics[width=.8\columnwidth]{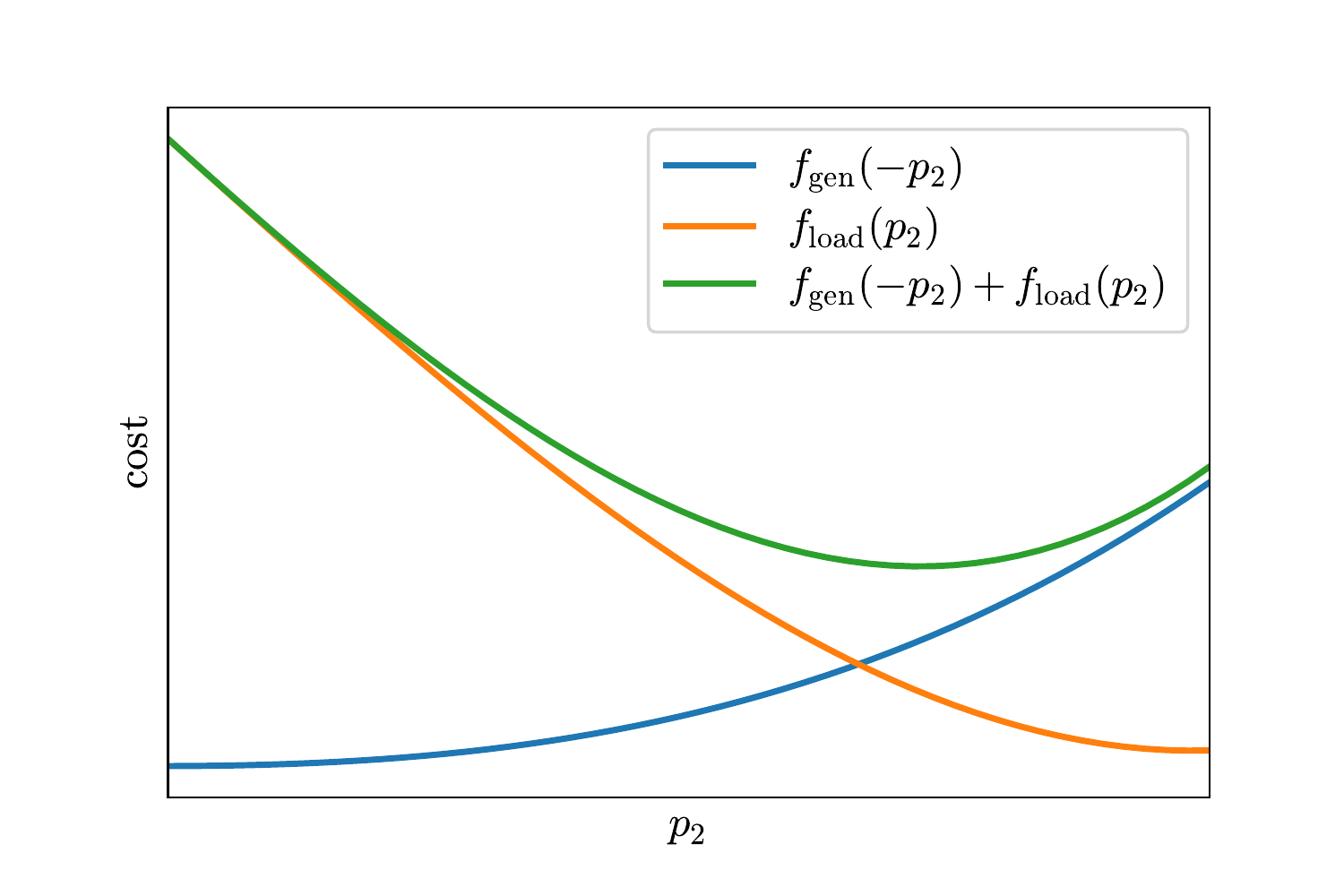}
  \caption{Generator, load, and total cost functions
  for the example of \S\ref{sec-ex-two-devices}.}
  \label{f-gen-costs}
\end{center}
\end{figure}

We can express this problem in a more natural form by
eliminating $p_1$ using $p_2 = -p_1$.  The power $p_2$ would typically be
positive, and corresponds to the power consumed by the load, which is the same 
as the power produced by the generator, $-p_1$.
The OPF problem is then to choose $p_2$ to minimize 
$f_{\rm load}(p_2) + f_{\rm gen}(-p_2)$.  
This is illustrated in figure~\ref{f-gen-costs},
which shows these two functions and their sum.

\begin{figure}
\begin{center}
  \includegraphics[width=.8\columnwidth]{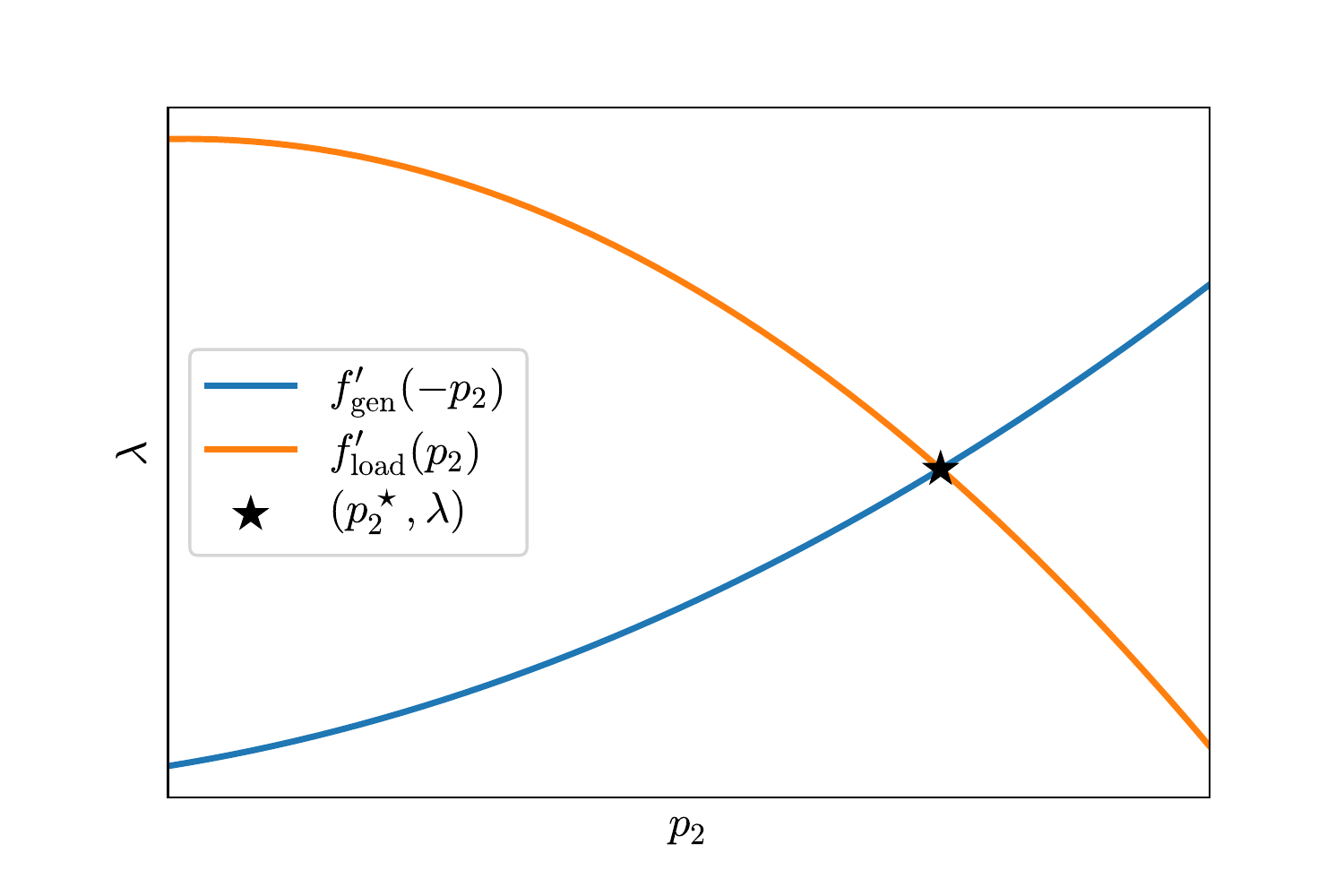}
  \caption{
  The derivatives of the generator and load cost functions
  are the supply and demand curves.
  The intersection of these curves is the point $(p_2^\star,\lambda)$.
  }
  \label{f-gen-load-sol}
\end{center}
\end{figure}

The optimality condition is $f_{\rm load}'(p_2) = f_{\rm gen}'(-p_2)= \lambda$.
We can interpret this as the crossing of the generator supply curve and the 
load demand curves.
This is shown graphically in figure~\ref{f-gen-load-sol},
with the optimal point given by the intersection
of the supply and demand functions at the point
$f_{\rm gen}'(p_2^\star) = f_{\rm load}'(-p_2^\star) = \lambda$.

\subsubsection{Three-bus example}
\label{s-three-bus-example}
We return to the three-bus example shown in
figure~\ref{f-transmission-network}. The units of power
are MW, with an implicit time period of one hour, and the 
units of payment are US dollars.

The two generators have quadratic cost functions.
The first has parameters $\alpha = 0.02$ $\rm \$/(\text{MW})^2$, 
$\beta = 30$ $\rm \$/(\text{MW})$,
$\pmin = 0$, and $\pmax = 1000$ $\text{MW}$.
The second has parameters $\alpha = 0.2$ $\rm \$/(\text{MW})^2$, 
$\beta = 0$ $\rm \$/(\text{MW})$,
$\pmin = 0$, and $\pmax = 100$ $\rm MW$.
Both of the loads are fixed loads;
the first consumes $\pfix = 50$ $\rm MW$,
and the second $\pfix = 100$ $\rm MW$.
The three transmission lines are lossless,
and have transmission limits
$\pmax = 50$ $\rm MW$, $\pmax = 10$ $\rm MW$, and $\pmax = 50$ $\rm MW$, 
respectively.

\begin{figure}
\begin{center}
  \scalebox{.8}{\begin{tikzpicture}

\def\termsep{2cm}

\newcommand{\priceNetOne}{\$ 33.6}
\newcommand{\priceNetTwo}{\$ 199.6}
\newcommand{\priceNetThree}{\$ 24.0}

\newcommand{\flowLoadOne}{50}
\newcommand{\flowLoadTwo}{100}

\newcommand{\payLoadOne}{}
\newcommand{\payLoadTwo}{ }

\newcommand{\flowGenOne}{-90}
\newcommand{\flowGenTwo}{-60}

\newcommand{\payGenOne}{ }
\newcommand{\payGenTwo}{ }

\newcommand{\flowLineOne}{50}
\newcommand{\flowLineOneBis}{-50}
\newcommand{\payLineOne}{ }

\newcommand{\flowLineTwo}{-10}
\newcommand{\flowLineTwoBis}{10}
\newcommand{\payLineTwo}{ }

\newcommand{\flowLineThree}{-50}
\newcommand{\flowLineThreeBis}{50}
\newcommand{\payLineThree}{ }

\node (n1)  [net, label=above:{\textcolor{price}{$\priceNetOne$}}] {\,net 1\,};

\node (t1)   [device, position=240:{\termsep} from n1, label=left:{\textcolor{pmnt}{$\payLineOne$}}] {line 1};
\node (t2)   [device, position=300:{\termsep} from n1, label=right:{\textcolor{pmnt}{$\payLineTwo$}}] {line 2};

\node (n2)  [net, position=240:{\termsep} from t1,
             label=left:{\textcolor{price}{$\priceNetTwo$}}] {\,net 2\,};
\node (n3)  [net, position=300:{\termsep} from t2,
             label=right:{\textcolor{price}{$\priceNetThree$}}] {\,net 3\,};

\node (t3)   [device, distance={\termsep}, right of=n2, label=below:{\textcolor{pmnt}{$\payLineThree$}}] {line 3};

\node (g1)   [device, left of=n1] {gen.\ 1};
\node (g2)   [device, below of=n3] {gen.\ 2};
\node (l1) [device, right of=n1] {load 1};
\node (l2) [device, below of=n2] {load 2};

\path[terminal] (n1) edge[pos=.50, above] node {\textcolor{flow}{$\flowGenOne$}}  (g1);
\path[terminal] (n1) edge[pos=.50, below] node {\textcolor{pmnt}{$\payGenOne$}}  (g1);

\path[terminal] (n1) edge[pos=.50, above] node {\textcolor{flow}{$\flowLoadOne$}}  (l1);
\path[terminal] (n1) edge[pos=.50, below] node {\textcolor{pmnt}{$\payLoadOne$}}  (l1);

\path[terminal] (n1) edge[pos=.50, left ] node {\textcolor{flow}{$\flowLineOne$}}  (t1);
\path[terminal] (n1) edge[pos=.50, left ] node {\textcolor{flow}{$\flowLineTwo$}}  (t2);
\path[terminal] (n2) edge[pos=.50, left ] node {\textcolor{flow}{$\flowLineOneBis$}}  (t1);
\path[terminal] (n3) edge[pos=.50, left ] node {\textcolor{flow}{$\flowLineTwoBis$}}  (t2);
\path[terminal] (n2) edge[pos=.50, below] node {\textcolor{flow}{$\flowLineThree$}}  (t3);
\path[terminal] (n3) edge[pos=.50, below] node {\textcolor{flow}{$\flowLineThreeBis$}}  (t3);

\path[terminal] (n2) edge[pos=.50, right] node {\textcolor{flow}{$\flowLoadTwo$}}  (l2);
\path[terminal] (n2) edge[pos=.50, left] node {\textcolor{pmnt}{$\payLoadTwo$}}  (l2);

\path[terminal] (n3) edge[pos=.50, left ] node {\textcolor{flow}{$\flowGenTwo$}}  (g2);
\path[terminal] (n3) edge[pos=.50, right ] node {\textcolor{pmnt}{$\payGenTwo$}}  (g2);

\end{tikzpicture}}
  \caption{The three-bus example, with solution;
  locational marginal prices in green
  and power flows in yellow.
    }
  \label{f-three-bus-sol}
\end{center}
\end{figure}
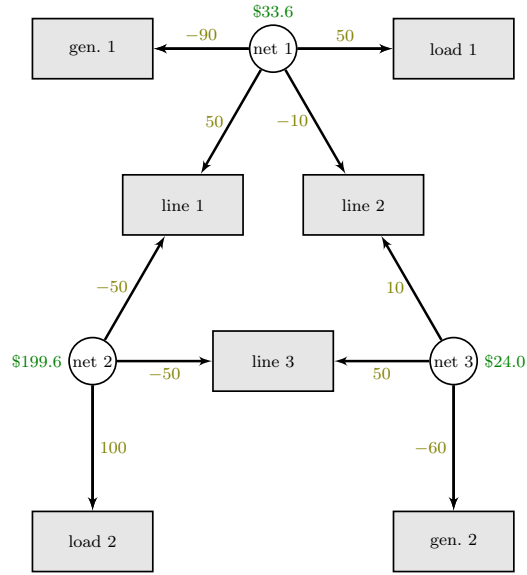

\paragraph{\emph{Results.}}
The solution of the OPF problem is shown in
figure~\ref{f-three-bus-sol},
and the payments to each device
are given in table 
\ref{t-three-bus-payments}.
The yellow numbers,
displayed next to the terminals,
show the optimal power flows $p^\star$.
The green numbers, next to the nets,
show the optimal prices $\lambda$.
First note that
the power flows into each net sum to zero,
indicating that this power flow satisfies
the conservation of power, \ie, condition (\ref{e-power-conservation}).
Furthermore,
all device constraints are satisfied:
the transmission lines transmit power according to their capacities,
each load receives its desired power,
and the generators supply positive power.

We see that power is cheapest near the second generator,
which is not near any load,
and expensive near the second load,
which is not near any generator.
Also note that although generator 2 produces
power more cheaply than generator 1,
the capacity limits of the transmission lines
limit production.
This has the effect that this generator is not paid much.
(See table~\ref{t-three-bus-payments}.)
Generator 1, on the other hand,
is paid much more,
which is justified by its advantageous proximity to load 1.
In addition to the generators,
the transmission lines earn payments by transporting power.
For example,
the third transmission line
is paid a substantial amount for transporting
power from net 3,
where power is generated cheaply,
to net 2, where there is a load, but no generation.
(The payment can be calculated
as the difference in price across its adjacent nets
multiplied by the power flow across it.)
So, we see that the two generators are paid for producing power,
the two loads pay for the power they consume, and the three 
transmission lines are paid for transporting power.
The payments balance out; they sum to zero.
The \texttt{cvxpower} code for this example is given in
appendix~\ref{s-cvxpower}.

\begin{table}
	\begin{center}
		\begin{tabular}{lr}
          Device     & Payment (\$) \\
          \hline
          \hline
          generator 1  & $          -  3024$ \\
          generator 2  & $          -  1440$ \\
          load 1   & $ \phantom{-} 1680$ \\
          load 2   & $  \phantom{-}19960$ \\
          line 1   & $          -  8300$ \\
          line 2   & $          -  96$ \\
          line 3   & $          -  8780$  \\
		\end{tabular}
	\end{center}
	\caption{The device payments
    for the three-bus example.}
	\label{t-three-bus-payments}
\end{table}

\section{Dynamic optimal power flow}
\label{c-dyn-opf}

\subsection{Dynamic network model}
\label{s-dyn-model}
In this section we generalize the static power flow model of \S\ref{c-stat-opf}
to dynamic optimization over $T$ time periods.
Each terminal has not just one power flow (as in the static case),
but instead a \emph{power schedule},
which is a collection of $T$ power flows,
each corresponding to one of the $T$ time periods.
Each device has a single cost function,
which associates a scalar cost with
the power schedules of its terminals.

The nets losslessly exchange power at each time period.
Power conservation holds if
the powers flowing into the net from the terminal devices sum to zero
for each of the $T$ time periods.
If this condition holds,
and if the cost associated with the terminal powers is finite,
then we say the powers are \emph{feasible}.

\subsubsection{Notation}
Much of the notation from the static case is retained.
However, in the dynamic case, the power flows are now described by a
matrix $p\in\reals^{M\times T}$.
The $m$th row of this matrix describes
the power schedule of terminal $m$.
The $t$th column describes the powers of all terminals
corresponding to time period $t$,
\ie, it is a snapshot of the power flow of the system 
at time period $t$.
The matrix $p_d = B_d p\in\reals^{M_d\times T}$
contains the power schedules of device $d$'s terminals.
The cost function of device $d$ is
$f_d:\reals^{M_d\times T}\to\reals\cup\{\infty\}$.
The system cost is again the sum of the device costs, \ie,
$f(p) = \sum_{d=1}^D f_d(p_d)$.
Conservation of power is written as the matrix
equation $Ap=0$ (\ie, $MT$ scalar equations).
Note that in the case of a single time period ($T=1$)
the dynamic case (and all associated notation) reduces to the static case.

\subsection{Optimal power flow}
\label{s-dyn-opf}
The \emph{dynamic optimal power flow problem} is
\begin{equation}
  \begin{array}{ll}
    \mbox{minimize} & f(p) \\
    \mbox{subject to} & Ap = 0,
  \end{array}
  \label{e-dyn-prob}
\end{equation}
with decision variable $p \in \reals^{M\times T}$,
which is the vector of power flows across all terminals at all time periods.
The problem is specified by the cost
functions $f_d$ of the $D$ devices,
the adjacency matrix $A$,
and the global-local matrices $B_d$, for $d=1,\ldots,D$.
Note that if $T=1$,
(\ref{e-dyn-prob}) reduces to the static OPF problem (\ref{e-static-prob}).

\subsubsection{Optimality conditions}
\label{s-dynamic-opt-conds}
If the cost function $f$ is convex and differentiable,
the prices are again the Lagrange multipliers of the
conservation of power constraint.
That is, the power flow matrix $p^\star$
and the price matrix $\lambda\in\reals^{N\times T}$
are optimal for (\ref{e-dyn-prob}) if and only if they satisfy
\[
\label{e-dyn-opt-cond}
\nabla f(p^\star) = A^T\lambda, \qquad
Ap^\star = 0.
\]
Note that $\nabla f(p)$ is a $M\times T$ matrix
of partial derivatives of $f$ with respect to the elements of $p$,
which means that the first equation consists of $MT$ scalar equations.
As in \S\ref{s-static-opt-conds},
these optimality conditions can be modified to handle the case of
a convex, nondifferentiable cost function (via subdifferentials),
or a nonconvex, differentiable cost function (in which case they 
are necessary conditions for local optimality).

\subsubsection{Solving the dynamic optimal power flow problem}
If all the device cost functions are convex,
then (\ref{e-dyn-prob}) is a convex optimization problem,
and can be solved exactly using standard algorithms,
even for large networks over many time periods
\cite[\S 1.3]{cvxbook}.
Without convexity of the device cost functions,
finding (and certifying)
a global solution to the problem is computationally intractable in general.
However, even in this case,
effective heuristics
(based on methods for convex optimization)
can often obtain practically useful (if suboptimal)
power flows, even for large problem instances.

\subsection{Prices and payments}
\label{s-dyn-lmps}
Here we explore the concept of \emph{dynamic} locational marginal prices.

\paragraph{\emph{Perturbed problem.}}
As in  \S\ref{s-static-LMP},
here we consider the perturbed system obtained by 
extracting a (small) additional amount of power into each net.
Now, however, the perturbation is a matrix
$\delta \in \reals^{N\times T}$,
\ie, the perturbation varies over the terminals and time periods.
The perturbed dynamic power flow problem is
\[
	\begin{array}{ll}
    \mbox{minimize} & f(p) \\
    \mbox{subject to} & Ap  + \delta = 0.
  \end{array}
\]
The optimal value of this problem,
as a function of $\delta$,
is denoted $F(\delta)$.

\paragraph{\emph{Prices.}}
Assuming the function $F$ is differentiable at the point $\delta = 0$,
the price matrix $\lambda $ is defined as
\[
\lambda = \nabla F(0)\in \reals^{N \times T}.
\]
This means that if a single unit of power is extracted 
into net $n$ at a time period $t$,
then the optimal value of problem (\ref{e-dyn-prob})
is expected to increase by approximately $\lambda_{nt}$.
Any element of $\lambda$ is, then, the  
marginal cost of power at a given network location 
over a given time period.
The $n$th row of this matrix is a vector of length $T$,
describing a \emph{price schedule}
at that net, over the $T$ time periods.
The $t$th column of this matrix is a vector of length $N$,
describing the prices in the system at all nets at time period $t$.
The price matrix $\lambda$ is also a Lagrange multiplier matrix:
together with the optimal power matrix $p^\star$,
it satisfies the optimality conditions (\ref{e-dyn-opt-cond}).

The price $\lambda_{nt}$
is the price of power at net $n$ in time period $t$.
If time period $t$ lasts for $h$ units of time,
then $\lambda_{nt}h$ is the price of energy during that time period.


\paragraph{\emph{Payments.}}
\label{s-dynamic-payment}
Here we extend the static payment scheme of \S\ref{s-static-payments}.
Device $d$ receives a total payment of
\[
\sum_{t=1}^T \lambda_{d,t}^T p_{d,t}
\]
over the $T$ time periods,
where $\lambda_d$ is the matrix of price schedules at nets
adjacent to device $d$,
\ie, $\lambda_d = B_dA^T \lambda$.
Note that we sum over the payments for each time period.
The sequence $(\lambda_{d,t}^T p_{d,t})$ for $t = 1, \ldots, T$,
is a \emph{payment schedule} or \emph{cash flow}.
In the dynamic case,
the payments clear at each time period and at each net,
\ie, all payments made at a single net sum to zero in each time period.
(This is a consequence of the fact that power is conserved at each time period,
at each net.)


\subsection{Profit maximization}
\label{s-dynamic-profits}
Under the payment scheme discussed above,
the profit of device $d$ is
\begin{equation*}
 -\sum_{t=1}^T \left( \lambda_{d,t}^T p_{d,t} \right) - f_d(p_d),
\end{equation*}
If $f_d$ is differentiable, this is maximized over $p_d$ if
\[
\nabla f_d(p_d) + \lambda_d = 0.
\]
Over all devices, this is precisely the first optimality condition of
(\ref{e-dyn-opt-cond}).
In other words,
the optimal power flow vector also maximizes
the individual device profits over the terminal power schedules of that device,
provided the adjacent net prices are fixed.
So, one can achieve network optimality
by maximizing the profit of each 
device, with some caveats.
(See, \eg, the recent work \cite{ma2018real}.)

\subsection{Dynamic device examples}
\label{s-dynamic-devices}
Here we list several examples of dynamic devices and their cost functions.
Whenever we list device constraints in the
device definition, we mean that the cost function is
infinite if the constraints are violated. 
(If we describe a device with only constraints,
we mean that its cost is zero if
the constraints are satisfied, and infinity otherwise.)

\subsubsection{Static devices}
All static devices,
such as the examples from \S\ref{s-static-devices},
can be generalized to dynamic devices.
Let $f_{d,t}(p_{d,t})$ be the (static)
cost of the device at time $t$.
Its dynamic cost is then the sum of all single-period costs
\[
f_d(p_d) = \sum_{t=1}^T f_{d,t}(p_{d,t}).
\]
In this case, we say that the device cost is \emph{separable} across time.
If all device costs have this property,
the dynamic OPF problem itself separates into 
$T$ static OPF problems; there are no constraints or objective terms
that couple the power flows in different time periods.
Static but time-varying devices can be used to model, for example,
time-varying availability of renewable power sources, or time-varying fixed loads.

\paragraph{\emph{Smoothing.}}
Perhaps the simplest generic example of a dynamic device objective that is not 
separable involves a cost term or constraint on the change in a power flow
from one time period to the next.   We refer to these generically as
\emph{smoothness penalties}, which are typically added to other, separable,
device cost functions.
A smoothness penalty has the form
\[
\sum_{t=0}^{T-1} \phi\left(p_{d,t+1} - p_{d,t}\right),
\]
where $\phi$ is a convex function.
(The initial value $p_{d,0}$ is a specified constant.)
Possible choices for $\phi$ are quadratic
($\phi(x) = \|x\|_2^2 = \sum_{i=1}^{M_d} x_i^2$),
absolute value or $\ell_1$
($\phi(x) = \|x\|_1 = \sum_{i=1}^{M_d} |x_i|$),
or an interval constraint function,
\[
	\phi_\text{ramp}(x) = 
\begin{cases}
	0 & -r_\text{down} \leq x \leq r_\text{up} \\
\infty & \mbox{otherwise,}
\end{cases}
\]
which enforces a maximum change in terminal power from one period to the next.
(These are called \emph{ramp rate limits} for a generator.)
For more details on smoothing, see \cite[\S 6.3.2]{cvxbook}.

\subsubsection{Dynamic generators}
\label{s-dynamic-generators}

\paragraph{\emph{Conventional generator.}}
An example of extending a static model
to a dynamic model is a conventional generator.
(See \S\ref{s-static-generators}.)
The cost function is extended to
\begin{equation}
\label{e-dynamic-controllable-gen}
f_d(p_d) =
\begin{cases}
	\sum_{t=1}^T \alpha p_{d,t}^2 + \beta p_{d,t}
 & p_{\rm min} \leq -p_{d,t} \leq p_{\rm max}, 
\quad t=1,\ldots,T  \\
\infty & \mbox{otherwise.}
\end{cases}
\end{equation}
where the scalar parameters $\alpha$, $\beta$, $p_\text{min}$,
and $p_\text{max}$ are the same as in  \S\ref{s-static-generators}.
(These model parameters could also vary over the $T$ time periods.)
We can also add a smoothing penalty or ramp rate limits, as discussed above.

\paragraph{\emph{Fixed generator.}}
Some generators cannot be controlled,
\ie, they produce a fixed power schedule $p_{\text{fix}}\in\reals^T$.
The device constraint is $-p_d = p_{\rm fix}$, $t=1, \ldots, T$.

\paragraph{\emph{Renewable generator.}}
Renewable generators can be controlled, with their
maximum power output depending on the availability
of sun, wind, or other power source.
That is, at each time period $t$,
a renewable generator can produce up to $p_{\text{avail}, t} \geq 0$
units of power.
The device constraint is that $-p_{d,t} \leq p_{{\rm avail}, t}$.

\subsubsection{Dynamic loads}

\paragraph{\emph{Fixed load.}}
A traditional load that consumes a fixed power in each period 
can be described by the simple constraints
\[
p_{d,t} = (p_\text{fix})_t,  \quad t=1, \ldots, T,
\]
where $p_\text{fix}$ is a fixed power schedule of length $T$.

\paragraph{\emph{Deferrable load.}}
A deferrable load requires a certain amount of energy over a given time window,
but is flexible about when that energy is delivered.
As an example,
an electric vehicle must be charged
before the next scheduled use of the vehicle,
but the charging schedule is flexible.
If the required energy over the time interval is $E_{\rm def} \geq 0$,
then the deferrable load satisfies
\[
\sum_{t=s}^e h p_{d,t} = E_{\rm def},
\]
where time periods $s$ and $e$ delimit the
start and end periods of the time window in which the device can use power,
and $h$ is the time elapsed between time periods.
We also require
\[
0 \leq p_{d,t} \leq p_\text{max}
\]
for time periods $t = s,\ldots, e$,
where $p_\text{max}$ is the maximum power the device can accept.
In addition,
we have $p_{d,t} = 0$ for time periods $t = 1,\dots, s-1$ and $t=e+1,\dots,T$.

\paragraph{\emph{Thermal load.}}
Here we model a temperature control system,
such as the HVAC
(heating, ventilation, and air conditioning)
system of a building,
the cooling system of a server farm,
or an industrial refrigerator.
The system has temperature $\theta_t$ at time period $t$,
and heat capacity $c$.
The system exchanges heat with its environment,
which has constant temperature $\theta_{\rm amb}$,
 and infinite heat capacity.
The thermal conductivity between the system and the environment is $\mu>0$.

We first consider the case of a cooling system,
such as a refrigerator or air conditioner,
with (cooling) coefficient of performance $\eta$.
The power used by the system at time $t$ is $p_{d,t}$.
The temperature changes with time according to
\[
\theta_{t+1} = \theta_{t}
+ (\mu/c) (\theta_{\rm amb} - \theta_t) - (\eta/c) p_{d,t}.
\]
The second term is the heat flow to or from the environment, and the third term
is heat flow from our cooling unit.

The temperature must be kept in some fixed range
\[
\theta_{\rm min} \leq \theta_t \leq \theta_{\rm max},
\quad t=1,\dots,T.
\]
The power consumption must satisfy the limits
\[
0 \leq p_{d,t} \leq p_{\rm max},
\quad t=1,\dots,T.
\]

The above model can be used to describe a heating system
when $\eta<0$.
In particular, for an electric heater,
$-\eta$ is the efficiency,
and is between $0$ and $1$.
For a heat pump, $-\eta$ is the (heating) coefficient of performance,
and typically exceeds $1$.
Other possible extensions include 
time-varying ambient temperature,
and higher-order models of the thermal dynamics.

\subsubsection{Storage devices}
We consider single-terminal storage devices,
including batteries, supercapacitors, flywheels, pneumatic storage,
or pumped hydroelectric storage.
We do not consider the specific details of these technologies,
but instead develop a simple model that
approximates their main features.

\paragraph{\emph{Ideal energy storage device.}}
We first model an ideal energy storage device, and then
specialize to more complicated models.
Let $E_t \in\reals_+$ be the  internal energy of the device
at the end of time period $t$.
(This is an internal device variable, whose value is fully specified
by the power schedule $p_d$).
The internal energy satisfies
\[
E_{t+1} = (1-\alpha) E_t + h p_{d,t}, \quad t=0,\dots,T-1,
\]
where $\alpha$ is the (per-period) \emph{leakage} rate, with
$0 < \alpha \leq 1$,
and $h$ is elapsed time between time periods.
The energy at the beginning of the first time period
is $E_0 = E_{\rm init}$,
where $E_{\rm init}$ is given.
We have minimum and maximum energy constraints
\[
E_{\rm min} \leq E_t \leq E_{\rm max}, \quad t=1,\dots, T,
\]
where $E_{\rm min}$ and $E_{\rm max}$ are the minimum 
and maximum energy limits.
In addition, we have limits on the charge and discharge
rate:
\[
p_\text{min} \leq p_d \leq p_\text{max}.
\]
We can impose constraints on the energy in a storage device, for example,
$E_T=E_{\rm max}$, \ie, that it be full at the last period.

The profit maximization principle allows us to relate prices in different time periods 
for an ideal lossless storage device (\ie, with $\alpha = 0$) that does not hit 
its upper or lower energy limit.  The prices in different periods at the adjacent 
net must be the same. (If not, the storage device could increase its profit by
charging a bit in a period when the price is low, and discharging the same energy when
the price is high.  This is analogous to a lossless transmission line that is not 
operating at its limits, which enforces equal prices at its two nets.   The transmission
line levels prices at two nets; the storage device levels the prices in two time 
periods.

\paragraph{\emph{Charge/discharge cost.}}
Many storage devices, such as batteries, degrade with use.
A charge/discharge usage penalty can be added to avoid overusing the device
(or to model amortization, or maintenance, cost).
We propose, in addition to the constraints above,
the cost
\begin{equation}
  \label{e-cost-ideal-storage}
  \beta \sum_{t=1}^T |p_{d,t}|,
\end{equation}
where $\beta$ is a positive constant,
and we are implicitly treating the vector $E$
as a function of the vector $p_d$.
For a battery whose capital cost is $C$,
and with an estimated
lifetime of $n_{\rm cyc}$ (charge and discharge) cycles,
a reasonable choice for $\beta$ is
\[
\beta= \frac{C h}{2 n_{\rm cyc}(E_{\rm max} - E_{\rm min})}.
\]

\begin{figure}
\begin{center}
  \scalebox{.8}{\begin{tikzpicture}

\node (net1)  [net] {};
\node (battery)  [device, position=180:{1.5cm} from net1] {ideal stor.};
\node (converter)  [device, position=0:{1.5cm} from net1] {lossy trans.};
\node (dX)  [position=000:{2cm} from converter] {\textcolor{flow}{}};
\node (label)  at ($(net1.north)+(-0.0cm,+0.8cm)$) {lossy battery};
\node (net2)  [position=0:{5.5cm} from net1] {};

\node (rightend) [position=000:{3cm} from net1] {};
\node (leftend)  [position=180:{3cm} from net1] {};

\path[terminal] (net1) edge[pos=.50, above] node
          {\textcolor{flow}{}}  (battery);
\path[terminal] (net1) edge[pos=.50, above] node
          {\textcolor{flow}{}}  (converter);
\path[terminal] (net2) edge[] node {}  (converter);

\draw[dashed] ($(battery.north west)+(-0.3cm,0.3cm)$)  rectangle
              ($(converter.south east)+(0.3cm,-0.3cm)$);

\end{tikzpicture}}
  \caption{
  A battery with conversion losses
  can be modeled as a compound device 
  made of an ideal storage unit and a lossy
  transmission line.
  }
  \label{f-lossy-battery}
\end{center}
\end{figure}
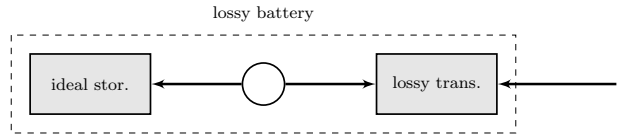

\paragraph{\emph{Conversion losses.}}
\label{s-lossy-battery}
In many cases, energy is lost when charging or discharging.
This can be modeled by adding
a lossy transmission line. 
(See \S\ref{par-lossy-transmission-line}.)
between the ideal storage device and its net,
as shown in figure~\ref{f-lossy-battery}.


%
%
%

\subsection{Home energy example}
\label{s-dyn-home-ener-example}
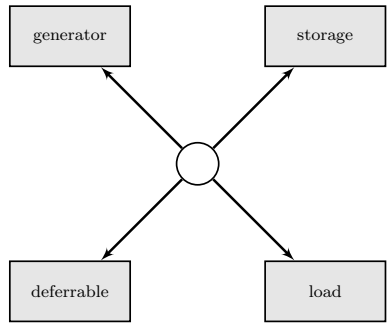
\begin{figure}
\begin{center}
\scalebox{.8}{\begin{tikzpicture}

\node (net)  [net] {};
\node (generator)  [device, above left of=net] {generator};
\node (storage)  [device, above right of=net] {storage};
\node (deferrable)  [device, below left of=net] {deferrable};
\node (load) [device, below right of=net] {load};

\node (rightend) [position=000:{3cm} from net] {};
\node (leftend)  [position=180:{3cm} from net] {};

\path [terminal] (net) -- (generator);
\path [terminal] (net) -- (storage);
\path [terminal] (net) -- (deferrable);
\path [terminal] (net) -- (load);

\end{tikzpicture}}
\caption{Network topology for the home energy example.}
\label{f-dynamic-example}
\end{center}
\end{figure}

We consider a home network with four devices:
a conventional generator, a fixed load, a deferrable load,
and an energy storage device.
They are all connected to a single net,
as shown in figure~\ref{f-dynamic-example}.
We operate them over a whole day, split into $T=1280$
time periods of 1 minute each.

\begin{figure}
\begin{center}
\includegraphics[width=.9\columnwidth]{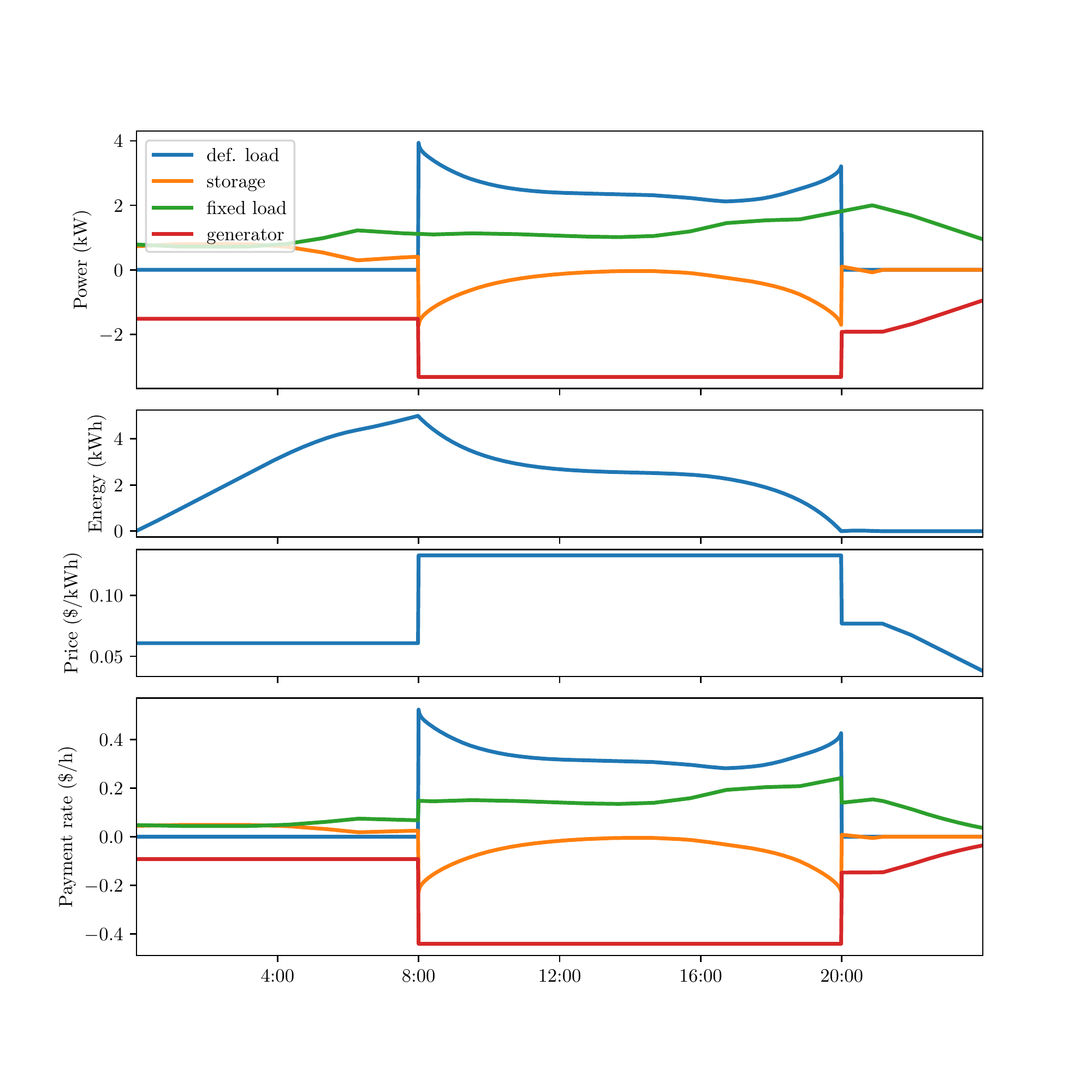}
\caption{
\emph{Top.}
Power consumption of the four devices
(with negative corresponding to production).
\emph{Upper middle.} Stored energy of the battery.
\emph{Lower middle.} Price at the net.
\emph{Bottom.} Device payment rate over time.
}
\label{f-home-energy-results}
\end{center}
\end{figure}

The generator cost function is given by
(\ref{e-quadratic-gen}),
with $\alpha = 0.0003$ $\rm \$/(kW)^2$, 
$\beta= 0$ $\rm \$/kW$,
$p_{\rm min} = 0$ $\rm kW$, and $p_{\rm max} = 6$ $\rm kW$.
The (ideal) energy storage device has
discharge and charge rates of
$p_{\rm min} =-2$ $\rm kW$, $p_\text{max} =  2$ $\rm kW$,
and minimum and maximum capacities 
$E_\text{min} = 0$ and
$E_\text{max} = 5$ $\rm kWh$, and is initially uncharged.
The deferrable load has a maximum power $p_{\rm max} = 5$ $\rm kW$,
and must receive $E_\text{def} = 30$ $\rm kWh$
of energy between 8:00 and 20:00.
The uncontrollable load has a time-varying power demand profile,
which is shown,
along with the problem solution,
in figure~\ref{f-home-energy-results}.

\paragraph{\emph{Results.}}
The optimal power and price schedules are shown in
figure~\ref{f-home-energy-results},
along with the internal energy of the storage device.
We see that the storage device and the deferrable load smooth out the
power demand over the time horizon,
and thus the total generation cost is reduced
(compared with the same system without 
a storage device, or with another fixed load instead
of the deferrable one).
The storage device charges (\ie, takes in energy) during the
initial time periods,
when power is cheap (because there is less demand),
and discharges (\ie, returns the energy) later, 
when it is more expensive.

When the deferrable load becomes active in time period 450 (\ie, 8:00),
there is even more flexibility in scheduling power,
and the price stays constant.
This is to be expected;
due to the quadratic cost of the generator,
the most efficient generation occurs
when the generator power schedule is constant,
and this can only happen if the price is constant.

The four payments by the devices are shown in table~\ref{t-dyn-payments}.  
We see that the generator is paid
for producing power, the fixed load pays for the power it consumes, 
as does the deferrable load (which pays less than if it were a fixed load).
The storage device is paid for its service, which is to transport power
across time from one period to another (just as a transmission line 
transports power in one time period from one net to another).

\begin{table}
	\begin{center}
		\begin{tabular}{lr}
          Device     & Payment (\$) \\
          \hline
          \hline
          generator  & $          -6.46$ \\
          deferable load  & $\phantom{-}3.99$ \\
          fixed load & $\phantom{-}2.83$ \\
          storage    & $          -0.36$ \\
		\end{tabular}
	\end{center}
	\caption{The device payments, summed over time.}
	\label{t-dyn-payments}
\end{table}

%
%
%
%

\section{Model predictive control}
\label{c-mpc}
The dynamic optimal power flow problem is useful for planning a schedule of power
flows when all relevant quantities are known in advance,
or can be predicted with high accuracy.
In many practical cases this assumption does not hold.  For example,
while we can predict or forecast future loads, the forecasts will not be 
perfect.  Renewable power availability is even harder to forecast.

\subsection{Model predictive control}
\label{s-mpc-def}
In this section we describe a standard method,
\emph{model predictive control} (MPC), also known as \emph{receding horizon control}
(RHC), that can be used to develop a real-time strategy or policy that
chooses a good power flow in each time period, and tolerates forecast errors.
MPC leverages our ability to solve a dynamic power flow problem over a horizon.
MPC has been used successfully in a very wide range of applications, including for the management of energy devices 
\cite{long2019generalised}.

MPC is a feedback control technique that naturally
incorporates optimization (\cite{bemporad2006model, mattingley2011receding}).
The simplest version is \emph{certainty-equivalent} MPC
with finite horizon length $T$, which is described below.
In each time period $t$, we consider a horizon that extends some fixed number of 
periods into the future, $t, t+1, \ldots, t+T-1$.
The number $T$ is referred to as the \emph{horizon} or \emph{planning 
horizon} of the MPC policy.
The device cost functions depend on various future quantities; while we know
these quantities in the current period $t$, we do not know them for future periods 
$t+1,t+2, \ldots, t+T-1$.
We replace these unknown quantities with predictions or forecasts,
and solve the associated dynamic power flow problem to produce a
(tentative) \emph{power flow plan}, that extends from the current
time period $t$ to the end of our horizon, $t+T-1$.
In certainty-equivalent MPC, the power flow plan is based 
on the forecasts of future quantities.
We then execute the first power flow in the plan, \ie, the power
flows corresponding to time period $t$ in our plan.
At the next time step, we repeat this process,
incorporating any new information into our forecasts.

To use MPC,
we repeat the following three steps at each time step $t$:
\begin{enumerate}
\item \emph{Forecast.}
  Make forecasts of unknown quantities to form an estimate of the
  device cost functions
  for time periods $t+1$, $t+2$, \dots, $t+T-1$.
\item \emph{Optimize.}
  Solve the dynamic optimal power flow problem (\ref{e-dyn-prob})
  to obtain a power flow plan
  for time periods $t, t+1, \dots, t+T-1$.
\item \emph{Execute.}
  Implement the first power flow in this plan,
  corresponding to time period $t$.
\end{enumerate}
We then repeat this procedure,
incorporating new information,
at time $t+1$.
Note that these steps can be followed indefinitely; the MPC method is always
looking ahead, or planning, over a horizon that extends $T$ steps into the future.
This allows MPC to be used to control power networks that run continuously.
We now describe these three steps in more detail.

\paragraph{\emph{Forecast.}}
At time period $t$,
we predict any unknown quantities relevant to system operation,
such as uncertain demand or the availability of renewable generators,
allowing us to form an approximate model of the
system over the next $T$ time periods.
These predictions can be crude,
for example as simple as a constant value such as the historical mean
or median value of the quantity.
The forecasts can also be sophisticated predictions based on previous values, historical
data, and even other, related quantities such as weather,
economic predictions, or futures prices (which are themselves market-based forecasts).
Appendix \S\ref{c-appendix-forecasts} describes a method
for creating basic forecasts, that are often adequate for MPC for dynamic energy 
management.

From predictions of these unknown quantities,
predictions of the device cost functions are formed
for time periods $t+1, \dots,t+T-1$.
At time $t$, we denote the predicted cost function
for device $d$ as $\hat f_{d|t}$.
The cost function for the entire system
is the sum of these cost functions, which we denote $\hat f_{|t}$.
(The hat above $f$ is a traditonal marker,
signifying that the quantity is
an estimate.)

\paragraph{\emph{Optimize.}}
We would like to plan out the power flows for the system
for time periods $t$ to $t+T-1$.
We denote by $p_{|t}$ the matrix of power flows
for all of the $D$ devices,
and for all of the $T$ time periods,
from $t$ to $t+T-1$.
We denote by $p_{\tau|t}$
the planned power flows for time period $\tau$.

To determine the planned power flows $p_{|t}$,
we solve the dynamic optimal power flow problem (\ref{e-dyn-prob}).
Using the notation of this section,
this problem is
\begin{equation}
  \begin{array}{ll}
    \mbox{minimize} & \hat f_{|t}(p_{|t}) \\
    \mbox{subject to} & Ap_{|t} = 0.
  \end{array}
  \label{e-mpc-prob}
\end{equation}
The variable is the planned power flow matrix $p_{|t}\in\reals^{M\times T}$.
The first column contains the power flows for the current period;
the second through last columns contain the planned power flows, based 
on information available at period $t$.

The optimization problem~(\ref{e-mpc-prob}) is sometimes augmented 
with \emph{terminal constraints} or \emph{terminal costs}, especially for storage devices.
A terminal constraint for a storage device specifies its energy level 
at the end of the horizon;
typical constraints are that it should be half full, or equal to the current 
value.  (The latter constraint means that over the horizon, the net total power 
of the storage device is zero.)
Without a terminal constraint, the solution of~(\ref{e-mpc-prob}) will have 
zero energy stored at the end of the horizon (except in pathological cases), 
since any stored energy could have been used to reduce some generator power,
and thereby reduce the cost.
A terminal cost is similar to a terminal constraint, except that it assesses
a charge based on the terminal energy value.

\paragraph{\emph{Execute.}}
Here, the first step of the planned power flow schedule
is executed,
\ie, we implement $p_{t|t}$.
(This could be as part of a larger simulation,
or this could be directly on the physical system.)
Note that the planned power flows $p_{t|t+1},\ldots, p_{t+T-1|t}$
are not directly implemented.
They are included for planning purposes only; their purpose is only
to improve the choice of power flows in the first step.

\subsection{Prices and payments}
\label{s-mpc-opt-conds}
Because the dynamic OPF problem
(\ref{e-dyn-prob}) is the same as problem (\ref{e-mpc-prob}),
the optimality conditions of \S\ref{e-dyn-opt-cond}
and the perturbation analysis of \S\ref{s-dyn-lmps}
also apply to (\ref{e-mpc-prob}),
which allows us to extend the concept of prices to MPC.
In particular, we denote the prices corresponding to a solution
of (\ref{e-mpc-prob}) as $\lambda_{|t}\in\reals^{M\times T}$.
This matrix can be interpreted as the predicted prices
for time periods $t+1,\dots,t+T-1$,
with the prediction made at time $t$.  (The first column contains the true
prices at time $t$.)


\paragraph{\emph{Payments.}}
\label{s-mpc-payments}
We can extend the payment scheme developed in
\S\ref{s-dynamic-payment} to MPC.
To do this,
note that the payment scheme in \S\ref{s-dynamic-payment}
involves each device making a sequence of payments
over the $T$ time periods.
In the case of MPC,
only the first payment in this payment schedule should be carried out;
the others are interpreted as planned payments.
Just as the planned power flows $p_{\tau|t}$
for $\tau = t+1,\ldots,t+T-1$
are never implemented,
but instead provide a prediction of future power flows,
the planned payments are never made,
but only provide a prediction of future payments.



\paragraph{\emph{Profit maximization.}}
In \S\ref{s-dynamic-profits},
we saw that given the predicted cost functions and prices,
the optimal power flows
maximize the profits of each device independently.
(We recall that we obtain the
prediction of
prices over the planning horizon
as part of the solution of
the OPF problem.)
Because the dynamic OPF problem is solved in each step of MPC,
this interpretation extends to our case.
More specifically,
given all information available at time $t$,
and a prediction of the prices $\lambda_{|t}$,
the planned power flows $p_{|t}$ maximize
the profits of each device independently.
In other words, 
if the managers (or owners) of 
each device agree on the predictions,
they should also agree that the planned power flows are fair.

We can take this interpretation a step further.
Suppose that at time $t$,
device $d$ predicts its own cost function as $f_{d|t}$,
and thus predicts the future prices to be $\lambda_{|t}$
(via the solution of the global OPF problem).
If the MPC of \S\ref{s-mpc-def} is carried out,
each device can be interpreted as carrying out MPC
to plan out its own terminal power flows to maximize its profit,
using the predicted prices $\lambda_{|t}$ during time period $t$.

\subsection{Wind farm example}
\label{s-mpc-wind-farm-example}
We consider a network consisting of a
wind generator,
a gas generator,
and a storage device,
and a fixed load,
all connected to one net.
The goal is to deliver a steady output of around $8$ $\rm MW$
to the fixed load,
which is the average of the available wind power over the month.
We consider the operation of this system for one month,
with each time period representing $15$ minutes.

\begin{figure*}
\begin{center}
\includegraphics[width=1\textwidth]{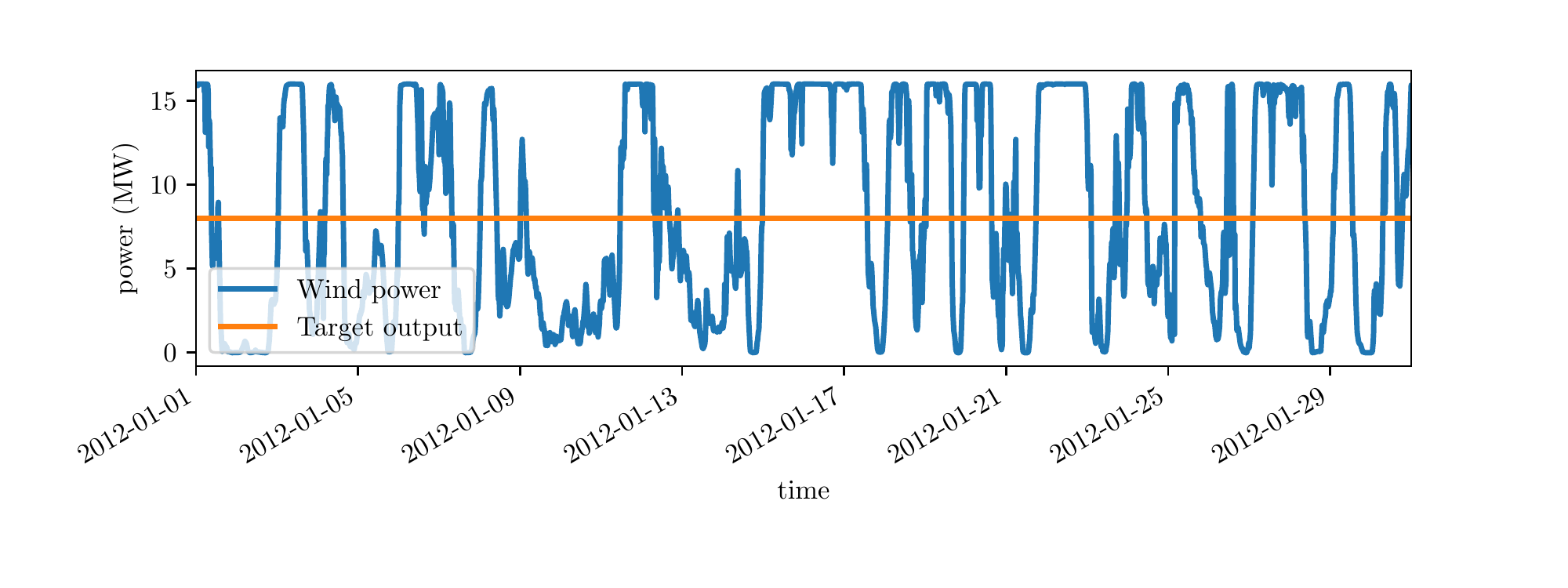}
\caption{
Available wind power for January 2012.
}
\label{f-wind-data}
\end{center}
\end{figure*}

The gas generator has the cost function given in \S\ref{s-dynamic-generators},
with parameters $\alpha=0.1$ $\rm \$ /(MW)^2$ and $\beta=20$ $\rm \$ /MW$.
The storage device has maximum charge and discharge rate of $5$ $\rm MW$,
and a maximum capacity of $50$ $\rm MWh$.
The wind generator is modeled as a renewable device,
as defined in \S\ref{s-dynamic-generators},
\ie, in each time period,
the power generated can be any nonnegative amount
up to the available wind power $p_{{\rm wind},t}$.
We show $p_{{\rm wind},t}$ as a function of the time period $t$
in figure~\ref{f-wind-data},
along with the desired output power.
The wind power availability data is provided by NREL (National Renewable 
Energy Laboratory), for a site in West Texas.
We solve the problem with two different methods, 
detailed below, and compare the results.
(Later, in \S\ref{s-robust-wind-farm}, we will introduce a third method.)

\paragraph{\emph{Whole-month dynamic OPF solution.}}
We first solve the problem as a dynamic OPF problem.
This requires solving a single problem that takes into account the entire month,
and also requires full knowledge of the available wind power.
This means our optimization scheme is \emph{prescient}, \ie, knows the 
future.
In practice this is not possible, but the performance in this prescient case
is a good benchmark to compare against, since no control scheme
could ever attain a lower cost.


\paragraph{\emph{MPC.}}
We then consider the practical case
in which the system planner does \emph{not}
know the available wind power in advance.
To forecast the available wind power,
we use the auto-regressive model developed in \S\ref{s-mpc-ar-model}, trained on data from the preceding year.
By comparing the performance of MPC
with the dynamic OPF simulation given above,
we get an idea of the value of (perfect) information,
which corresponds to the amount of additional cost incurred
due to our imperfect prediction of available wind power.

\begin{figure*}
\begin{center}
\vspace{-1cm}
\includegraphics[width=1.0\textwidth]{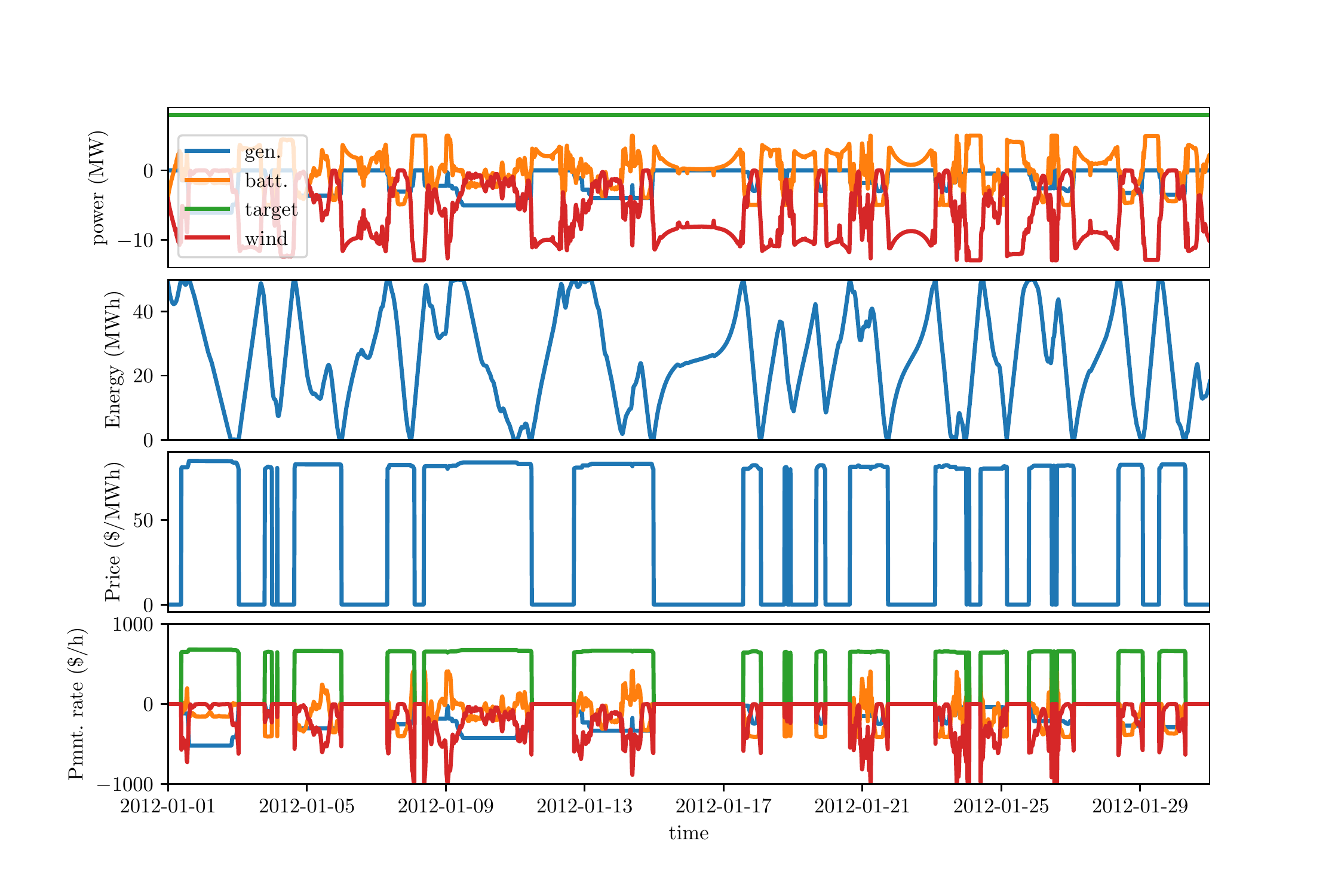}
\includegraphics[width=1.0\textwidth]{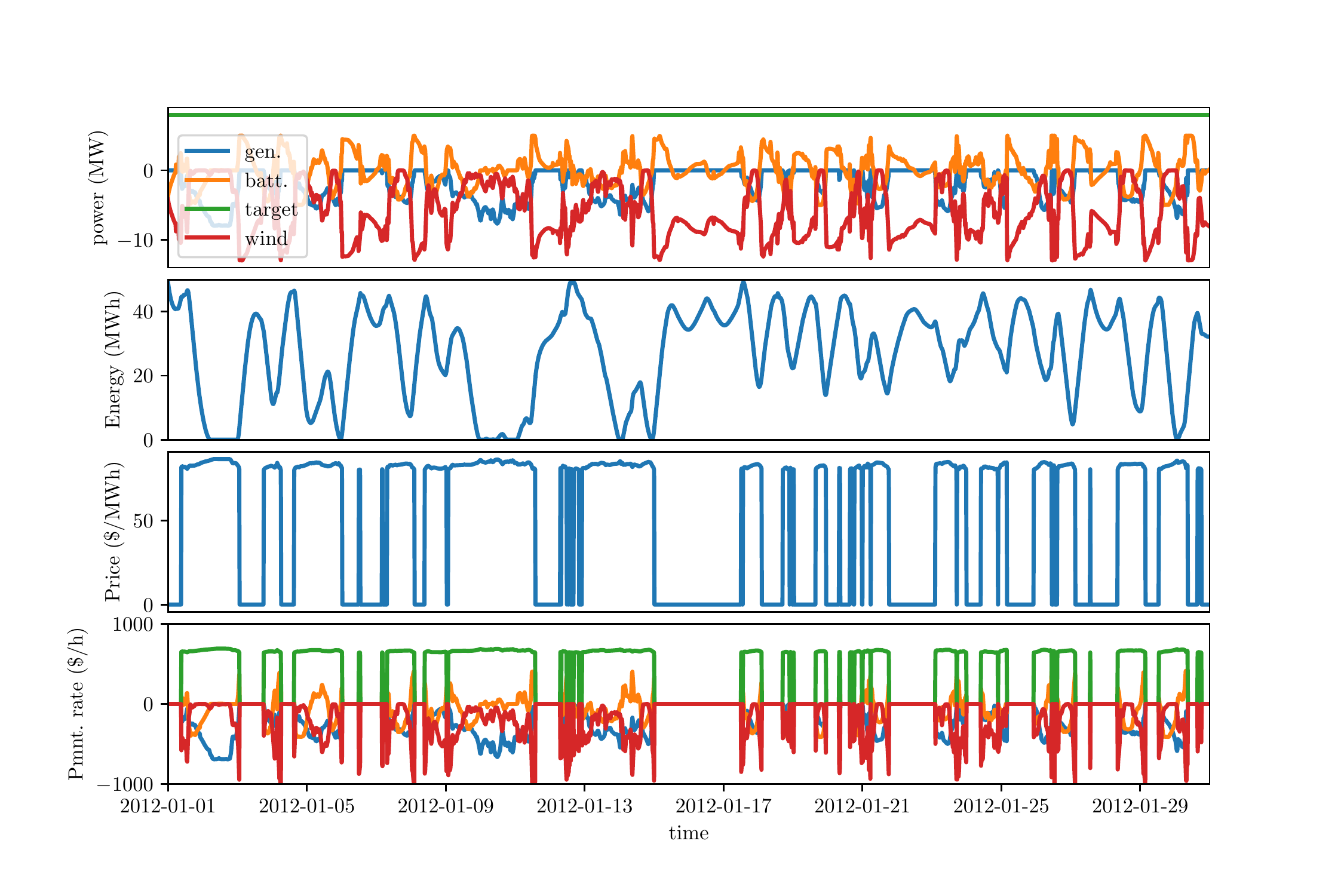}
\caption{
\emph{Top.} The full-month dynamic OPF solution.
\emph{Bottom.} The MPC simulation.
}
\label{f-wind-results}
\end{center}
\end{figure*}

\paragraph{\emph{Results.}}
The power flows obtained by solving the problem
using dynamic OPF and MPC are shown in figure~\ref{f-wind-results}.
The values of the cost function obtained using
dynamic OPF and MPC
were $\$3269$ and $\$3869$, respectively.
This difference reflects the cost of uncertainty,
\ie, the difference gives us an idea
of the value of having perfect predictions.
In this example the difference is not negligible, and suggests that investing in better 
wind power forecasting could yield greater efficiency.

In figure~\ref{f-wind-results},
we also show the prices (in time), as well as the payments made by each device.
Note that the price is ``set'' by the power production of the gas generator.
(This is because the price is given be the derivative of the cost functions
of adjacent devices.)
This means that when the gas generator produces power,
the price is positive;
otherwise, it is zero.
Also note that when the price is zero,
no payments are made.

In table~\ref{t-wind-payments},
we show the total payment of each device,
using dynamic OPF and MPC.
We see that under the dynamic OPF method,
the storage device is paid more than under MPC.
This is because storage is more useful precisely
when forecasting is accurate.
(For example, with no knowledge of future wind power availability,
the storage device would not be useful.)
Similarly, the gas generator is paid more under MPC.
This makes sense;
a dispatchable generator is more valuable if
there is more uncertainty about future renewable power availability.

\begin{table}
	\begin{center}
		\begin{tabular}{lrrr}
          Device     & Dynamic OPF & MPC  \\
          \hline
          \hline
          wind generator  & $ -98.4$ & $-115.3$ \\
          storage         & $ -54.0$ & $ -36.0$ \\
          load            & $ 255.6$ & $ 273.6$ \\
          gas generator   & $-103.0$ & $-122.3$ \\
		\end{tabular}
	\end{center}
	\caption{The total device payments over the entire month,
      in thousands of dollars.}
	\label{t-wind-payments}
\end{table}


\section{Optimal power flow under uncertainty}
\label{c-robust-mpc}
In this section,
we first extend the dynamic model of
\S\ref{c-dyn-opf} to handle uncertainty.
We do this by considering multiple plausible predictions of uncertain values,
and extending our optimization problem to handle multiple predictions or forecasts.
We will see that prices extend naturally to the uncertain case.
We then discuss how to use the uncertain optimal power flow problem
in the model predictive control framework of \S\ref{c-mpc}.

\subsection{Uncertainty model}
\label{s-uncert-model}

\paragraph{\emph{Scenarios.}}
Our uncertainty model considers $S$ discrete \emph{scenarios}.
Each scenario models a distinct contingency,
\ie, a different possible outcome
for the uncertain parameters in the 
network over the $T$ time periods.
The different scenarios can differ in the values at different time periods of 
fixed loads,
availability of renewable generators, and even the capacities of 
transmission lines or storage devices.
(For example, a failed transmission line has zero power flow.)

\paragraph{\emph{Scenario probabilities.}}
We assign a probability of realization to each scenario,
$\pi^{(s)}$ for $s=1,\dots,S$.
For example, we might model a \emph{nominal} scenario
with high probability, and a variety of \emph{fault}
scenarios, in which critical components break down, 
each with low probability.
The numbers $\pi^{(s)}$ form a probability distribution
over the scenarios,
\ie, $\pi^{(s)}\geq 0$ and $\sum_{s=1}^S \pi^{(s)}=1$.

\paragraph{\emph{Scenario power flows.}}
We model a different network power flow for each scenario.
The power flows for all terminals, time periods, and scenarios
form a (three-dimensional) array $p\in\reals^{M\times T\times S}$.
For each scenario $s$
there is a power flow matrix $p^{(s)}\in\reals^{M\times T}$,
which specifies the power flows on each of the $M$ terminals
at each of the $T$ time periods, under scenario $s$.
From the point of view of the system planner,
these constitute a \emph{power flow policy},
\ie, a complete contingency plan consisting of a power schedule
for each terminal, under every possible scenario.

We refer to the vector of powers for a device as
$p_d\in\reals^{M_d\times T\times S}$,
where $M_d$ is the number of terminals for device $d$.
This array can be viewed as a power flow policy specific to device $d$.
We denote by $p_d^{(s)}\in\reals^{M_d\times T}$
the submatrix of terminal power flows
incident on device $d$ under scenario $s$.

As before,
for each time period, and under each scenario,
the power flows incident on each net sum to zero:
\[
Ap^{(s)} = 0, \quad s=1,\dots,S.
\]
In the case of a single scenario ($S = 1$)
$p$ is a $M\times T$ matrix,
which corresponds to the power flow matrix
of \S\ref{c-dyn-opf}.

\paragraph{\emph{Scenario device cost functions.}}
The device cost functions may be different under each scenario.
More specifically, under scenario $s$,
device $d$ has cost $f_d^{(s)}$,
such that
$f_d^{(s)}:\reals^{M_d\times T}\to\reals\cup\{\infty\}$.
Note that the network topology,
including the number of terminals for each device,
does not depend on the scenario.
We define the cost function of device $d$ as its expected cost over all scenarios
\[
f_d(p_d) = \sum_{s=1}^S \pi^{(s)} f_d^{(s)}\big(p_d^{(s)}\big).
\]
In the case of a single scenario,
this definition of device cost coincides with the definition
given in \S\ref{c-dyn-opf}.
The expected total system cost is the sum of the expected device costs
\[
f(p) = \sum_{d=1}^D f_d(p_d).
\]

\subsection{Dynamic optimal power flow under uncertainty}

So far the different scenarios are not coupled, except possibly by
having common starting values for smoothness contraints.
To minimize $f$ subject to the power conservation constraint, we 
solve the $S$ dynamic OPF problems associated with each of the scenarios.

We are now going to couple the power flows for different 
scenarios with
an \emph{information pattern constraint}, which states that in time 
period $t=1$, the power flows for all $S$ scenarios must agree,
\ie, $p^{(1)}_{1} = \cdots = p^{(S)}_{1}$.
The \emph{uncertain dynamic optimal power flow problem} is
\begin{equation}
  \begin{array}{lll}
    \mbox{minimize} & f(p) \\
    \mbox{subject to} & Ap^{(s)} = 0,
                          &\quad s=1,\dots,S \\
                      & p^{(1)}_{1} = \cdots = p^{(S)}_{1},
  \end{array}
  \label{e-uncert-prob}
\end{equation}
where the variables are the scenario power flows
$p^{(s)}\in\reals^{M\times T}$, for $s=1,\dots,S$.
We can describe this problem as follows.   We create a full power flow
plan for each scenario, with the constraint that the first period power 
flow must be the same in each of the scenarios.

In terms of stochastic control, or optimization with recourse, the 
information pattern constraint corresponds to a 
very simple \emph{information pattern}, 
which is a description of what we know before we decide 
on our action.  
We have $S$ scenarios, one of which will occur; we must 
make the first choice, \ie, decide the current power flows \emph{before} we 
know which of the $S$ scenarios will actually occur. At period $t=2$,
the scenario that obtains is revealed to us.
Of course we do not believe this model, since the scenarios are just a (very small)
sampling of what might actually occur, and the future is not in fact revealed to us in
entirety at period $2$.   This is simply a heuristic for making good choices of 
the current period power flows that takes into account the fact that the future 
is uncertain.

\subsection{Prices and payments}
\label{s-uncert-LMP}
We now discuss locational marginal prices under our uncertainty model.
Suppose we inject extra power into each net, at each point in time,
for each scenario.
We describe these injections by scenario-specific matrices
 $(\delta^{(1)}, \ldots, \delta^{(S)})\in\reals^{N\times T\times S}$.
Power conservation, for each scenario, requires
\begin{equation}
Ap^{(s)} + \delta^{(s)} = 0, \quad s = 1, \ldots, S,
\label{e-uncert-perturbed-conservation}
\end{equation}
\ie, the extra power injected into each net,
summed with all power outflows along the incident terminals,
is zero.
If we solve problem (\ref{e-uncert-prob}),
with the power conservation constraints
replaced by the perturbed equations (\ref{e-uncert-perturbed-conservation}),
the optimal cost will change
to reflect the amount of power injected into each net;
we define $F(\delta^{(1)}, \ldots, \delta^{(S)})$ as 
the optimal value of the perturbed problem,
when the power injected under each scenario are 
$(\delta^{(1)}, \ldots, \delta^{(S)})$.
Note that $F(0)$ is the optimal value of the original, unperturbed problem.

Then, the price matrices $(\lambda^{(1)}, \ldots, \lambda^{(S)})\in\reals^{N\times T\times S}$, for each scenario, satisfy
\begin{equation}
(\pi^{(i)}\lambda^{(1)}, \ldots, \pi^{(S)}\lambda^{(S)}) = \nabla F(0).
\end{equation}
This means that the prices are given by the gradient $F$,
scaled up by the reciprocals of the scenario probabilities.
These matrices represent the predicted price of power at each net, 
each point in time, under each scenario.

It can be shown that the prices respect a constraint
similar to the information pattern constraint discussed above,
\ie, the prices can be chosen to coincide for the first time period,
across all scenarios, \ie,
\[
\lambda^{(1)}_{1} = \cdots = \lambda^{(S)}_{1},
\]
where $\lambda_{1}^{(s)}$ is the vector of prices during the first
time period, under scenario $s$. This property is important for the
payment scheme.

\paragraph{\emph{Payments.}}
We can extend our payment scheme from the dynamic case
to the uncertain case.
Each device's expected payment is
\[
\sum_{s=1}^S \pi^{(s)} \lambda_d^{(s)T} p_d^{(s)},
\]
where the expectation is taken over the various scenarios,
with the predicted price trajectories. 
 
Note that if we operate in the model predictive control framework of \S\ref{c-mpc},
only the first step payment is actually carried out.
In the next-step information pattern the first period prices coincide under each scenario,
so the payment does not the depend on the scenarios.

\paragraph{\emph{Profit maximization.}}
Given optimal power flows and prices,
the scenario power flows maximize the expected profit
for any device $d$
\[
\sum_{s=1}^S \pi^{(s)} 
\Big(\lambda^{(s)T} p_d^{(s)} - f_d^{(s)}\big(p_d^{(s)}\big) \Big),
\]
subject to an appropriate information pattern constraint.
This can be interpreted as follows.
Each device maximizes its own expected profit,
using the same uncertainty model (\ie, the scenario costs and 
 probabilities) as the system planner.
Note that each device maximizes its expected profit,
without caring about its variance, as is customary in 
model predictive control.
In the language of economics,
each device is assumed to be risk-neutral.
(One could include risk aversion in the cost function
of problem (\ref{e-uncert-prob}),
using a concave utility function;
see, \eg, \cite[\S11.5]{luenberger1995microeconomic}.)

\subsection{Robust model predictive control}
\label{s-robust-mpc-def}
Here we introduce an extension of the MPC framework
presented in \S\ref{c-mpc} to handle prediction uncertainty.
During the predict stage,
we allow for multiple forecasts of uncertain values.
We then plan the power flows by solving (\ref{e-uncert-prob}),
with each scenario corresponding to a forecast.
We repeat the following three steps at each time step $t$.
\begin{enumerate}
\item \emph{Predict.}
  We make $S$ plausible forecasts of unknown future quantities.
  Each forecast is a scenario,
  to which we assign a probability of occurrence.
  For each forecast, we form appropriate device cost functions.
\item \emph{Optimize.}
  We plan the power flows for each scenario
  by solving problem \eqref{e-uncert-prob},
  so the first planned power flows coincide under all scenarios.
\item \emph{Execute.}
  We execute the first power flow in this plan,
  \ie, the one corresponding to time period $t$
  (which coincide under all scenarios).
\end{enumerate}
We then repeat this procedure,
incorporating new information,
at time $t+1$.
We now describe these three steps in more detail.

\paragraph{\emph{Predict.}}
At time period $t$,
we make $S$ forecasts of all unknown quantities relevant to system operation.
Typically, these forecasts are generated using
a stochastic model of future variables.
(For example,
we can use a statistical model to generate several realistic generation
profiles for a solar generator, over the course of one day.)
We discuss some ideas for modeling in appendix~\ref{c-appendix-forecasts}.

Each forecast corresponds to a scenario;
for each forecast,
we form a scenario device cost function.
Under scenario $s$,
we denote the cost function for device $d$ as
$\hat f_{d|t}^{(s)}$.
The cost function for the entire system
is denoted $\hat f_{|t}^{(s)}$.

\paragraph{\emph{Optimize.}}
We  plan the system power flows 
for time periods $t$ to $t+T-1$, under each scenario, 
by solving the dynamic optimal power flow problem
with uncertainty
(\ref{e-uncert-prob}).
Denoting by $p_{|t}^{(s)}$ the matrix of power flows
under scenario $s$
for each of the $T$ future time periods $t+1, \ldots, t+T$, 
we solve
\begin{equation}
  \begin{array}{ll}
    \mbox{minimize} & \hat f(p_{|t}) \\
    \mbox{subject to} & Ap_{|t}^{(s)} = 0
                          \quad s=1,\dots,S \\
                      & p_{t +1|t}^{(s)} = p_{|t}^{\rm nom}
                          \quad s=1,\dots,S,
  \end{array}
  \label{e-uncert-prob-mpc}
\end{equation}
where the variables are the planned power flow matrices
$p_{|t}^{(s)}\in\reals^{M\times T}$, for each scenario $s$,
and the common first power flow, $p_{|t}^{\rm nom}$.

\paragraph{\emph{Execute.}}
The first step of the planned power flow schedule
is executed, \ie, we implement $p_{|t}^{\rm nom}$.
Note that the planned power flows
$p_{t+1|t}^{(s)},\ldots p_{t+T-1|t}^{(s)}$
are never directly implemented.
They are included for planning purposes,
under the assumption that
planning out $T$ steps, and under $S$ scenarios,
increases the quality of the first step in that plan.

\paragraph{\emph{Prices, payments, and profit maximization.}}
As noted in \S\ref{s-uncert-LMP},
the prices can be chosen such that
all prices for the first period coincide.
In the notation of MPC, we call these prices $\lambda_{|t}$.
These prices can be made the basis for a payment scheme.
At time period $t$, device $d$ is paid
\[
\lambda_{d|t}^T p_{d|t},
\]
where $\lambda_{d|t}$ is the vector of prices corresponding
to all nets adjacent to device $d$.
As in the static and dynamic case,
this payment scheme has the property
that the sum of all payments made is zero,
\ie, the payment scheme is revenue neutral.

In addition, the argument about 
profit maximization under (standard) MPC 
in \S\ref{s-mpc-payments} extends to the
robust MPC setting. 
If we assume that the managers of the
devices agree on the cost functions 
and probabilities of the different scenarios,
then they should agree that the planned power
flows are fair, and each device $d$ 
maximizes
its own expected profit (disregarding risk) 
by implementing the
optimal power $p_{d|t}^{nom}$.

\subsection{Uncertain device examples}

\paragraph{\emph{Deterministic device.}}
Any (deterministic) dynamic device
can be extended to an uncertain device.
Such a device has an identical cost function under each scenario.

\paragraph{\emph{Renewable generator.}}
Many renewable generators, such as solar and wind generators,
have uncertain energy production.
In this case,
the generator produces a potentially different amount of power
under each scenario.
If the generator produces power $P_t^{(s)}$ at time $t$
under scenario $s$,
then the generator power $p_{d,t}^{(s)}$ at time $t$
under scenario $s$ is
\[
0 \leq p_{d,t}^{(s)} \leq P_{t}^{(s)}, \quad t = 1,\dots,T.
\]

\paragraph{\emph{Uncertain load.}}
Loads can have an uncertain consumption pattern.
We assume an uncertain load consumes $P_t^{(s)}$
at time $t$ under scenario $s$.
This means that the power flows satisfy
\[
p_{d,t}^{(s)} = P_{t}^{(s)}, \quad t = 1,\dots,T, \quad s=1,\ldots, S.
\]

\paragraph{\emph{Unreliable transmission line.}}
Recall the definition of transmission line from \S\ref{c-stat-opf}.
Under all scenarios for which the transmission line works,
the device cost function is as described in \S\ref{s-transmission-lines}.
For scenarios in which the transmission line fails,
the two terminal power flows must both be zero,
\ie, we have $p_d^{(s)}=0$.

\subsection{Wind farm example}
\label{s-robust-wind-farm}
We extend the example of \S\ref{s-mpc-wind-farm-example} with uncertain
predictions of the wind power available. 
The network consists of a wind generator, a gas generator, 
and a storage device, all connected to one net. 
We consider the operation of this system for one month, 
with each time period representing 15 minutes.
The uncertain MPC example of \S\ref{s-mpc-wind-farm-example} uses
a single prediction of the future wind power available, obtained with
the AR model described in \S\ref{s-mpc-ar-model}.
To apply robust MPC, we require multiple forecasts of the unknown quantity.
We use the framework of \S\ref{c-appendix-forecasts}
to obtain $K=20$ such forecasts.

\begin{figure*}
\begin{center}
\includegraphics[width=.9\textwidth]{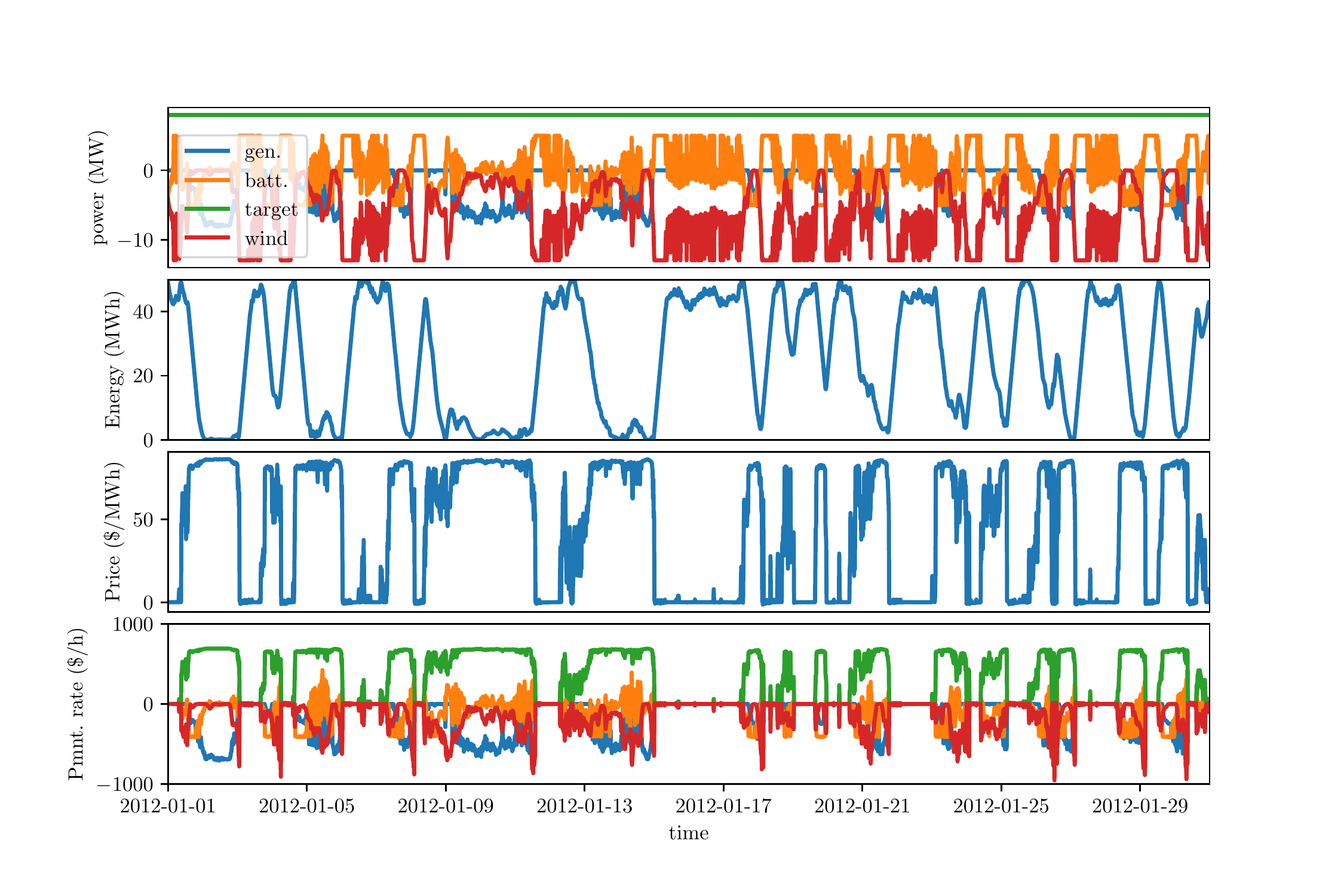}
\caption{
The robust MPC simulation.
}
\label{f-wind-rmpc}
\end{center}
\end{figure*}

\paragraph{\emph{Results.}}
The power flows obtained using robust MPC, with $K=20$ scenarios,
each with a different prediction of the uncertain wind power available
(and all other parameters equal), are shown in figure~\ref{f-wind-rmpc}.
The value of the cost function is \$$3291$.
This is not much higher than \$$3269$, the cost obtained using DOPF.
This illustrates a key point:
Even though our predictions of the wind power are fairly inaccurate,
the performance of the resulting control scheme is similar
to one that uses perfect predictions.

In table~\ref{t-robust-wind-payments},
we show the total payment of each device,
under the robust MPC formulation,
as well as the MPC formulation
discussed in \S\ref{s-mpc-wind-farm-example}.
The pattern here is similar to table~\ref{t-wind-payments};
the storage becomes more useful,
and is therefore paid more, when the forecasts are accurate,
and the gas generator is paid less.

\begin{table}
	\begin{center}
		\begin{tabular}{lrr}
          Device    &  MPC & Robust MPC \\
          \hline
          \hline
          wind generator  & $-115.3$ & $ -98.1$ \\
          storage         & $ -36.0$ & $ -48.9$ \\
          load            & $ 273.6$ & $ 251.6$ \\
          gas generator   & $-122.3$ & $-104.5$ \\
		\end{tabular}
	\end{center}
	\caption{The total device payments over the entire month,
      in thousands of dollars.}
	\label{t-robust-wind-payments}
\end{table}

\subsection*{Acknowledgments}
This research was partly supported by MISO energy;
we especially thank Alan Hoyt and DeWayne Johnsonbaugh of MISO for 
many useful discussions.



\bibliographystyle{plain}
\bibliography{dyn_ener_man}



\section{Appendix: Forecasts}
\label{c-appendix-forecasts}
Forecasts or predictions of time series, such as
loads or power availability of renewable generators,
are critical components of the MPC formulation
of \S\ref{c-mpc} or the robust MPC formulation of \S\ref{c-robust-mpc}.
We have already noted that the forecasts do not need to be very good
to enable MPC or robust MPC to yield good performance.
Even simple forecasts of quantities, such as predicting that future
values will simply be equal to the current value, can give reasonable
MPC performance in some cases.

Time series modeling and forecasting is a well studied problem 
in a variety of fields, such as statistics \cite{friedman2001elements}, 
machine learning, and econometrics \cite{shumway2006time}. 
These and many other references describe sophisticated forecasting techniques
that can be used.

In this section we describe a simple method to
forecast a scalar time series.  Our simple model takes into account
seasonal variations, both long and short term, in a baseline time
series, that depends only on time.
It also takes into account short-term deviations from the baseline,
based on recent past deviations.
Similar techniques can
be applied to vector time series, either by separating them
into their scalar components,
or by swapping the vector of model parameters
with an appropriate matrix of parameters.

\paragraph{\emph{General model fitting methods.}}
The simplest method to fit a model from data uses basic 
least squares or regression \cite{vmlsbook};
more sophisticated methods based on convex optimization use 
loss functions and regularizers
that give robust estimates, or sparse parameters (\ie, regressor selection)
\cite[Chap.~6]{cvxbook}.
Much more sophisticated forecasts can be developed using advanced techniques
like random forest models or neural networks \cite{friedman2001elements}.
We recommend starting with simple forecasts (even the constant one described above) 
and slowly increasing the complexity and sophistication of the forecast and 
evaluating the improvement (if any) on the control performance using MPC or 
robust MPC.
In a similar way we recommend starting with least squares fitting techniques before
moving to more sophisticated methods.


\subsection{The baseline-residual forecast}
We consider a time series $x_t \in \reals$, where 
the index $t=1, 2, \ldots$ represents time or period.
We think of $t$ as the current time; $t-1$ then refers to the previous period,
and $t+2$ refers to the period after the next period.
The series might represent the power availability of a renewable generator, 
or the power of a fixed load, with $t$ representing, \eg, the hours or 5 minute periods.
At any time $t$ we assume we have access to the current and
past observations 
 \[
 x_t, x_{t-1}, x_{t-2}, \ldots
 \]
and we are to make future predictions 
\[
\hat x_{t+1 | t}, \hat x_{t+2 | t}, \ldots, \hat x_{t+T-1| t},
\]
where by $\hat x_{t+\tau|t}$ we mean the prediction of the value at time
$t+\tau$, made at time $t$. (We use superscript $\hat{\ }$ to denote a predicted, and not actual, value.)
In the notation $t+\tau|t$, the lefthand part $t+\tau$ refers to the time
of the quantity being forecast; the righthand part $t$ refers to the time
at which the forecast is made.  We can interpret $\tau$ as the number of 
time periods into the future we forecast; for example, $\hat x_{t+1|t}$ is
called the \emph{one-step-ahead predictor} or forecast.
It is our prediction of the next value of the series, given what we know
now.

Using our forecast notation, we can express the simple constant forecast
as $\hat x_{t+\tau|t}=x_t$.  This predicts that all future values 
will be the same as the current value.  While this is rarely a good forecast 
for a time series (unless the time series is very slowly changing) it can
be adequate for MPC.
In the next few subsections below we describe forecasts that are 
a bit more sophisticated than the constant forecast, and often work very well.

\paragraph{\emph{Baseline-residual decomposition and forecast.}}
We model the time series as the sum of two components: 
a seasonal \emph{baseline} $b_t \in \reals$, which takes
into account variations due to, \eg, hourly, daily, annual, and weekly seasonalities 
and periodicities, 
and a \emph{residual} $r_t$, which is the deviation from the baseline,
\[
x_t = b_t + r_t, \quad t = 1, \ldots.
\]
The residual time series is also sometimes called the \emph{seasonally adjusted}
or \emph{baseline adjusted} time series.  It tells us how much larger or smaller
the values are, compared to the baseline.
We fit the baseline $b_t$ using some past or historical data,
as explained below.

Our predictions have the form
\[
\hat x_{t+\tau | t} = b_{t+\tau} + \hat r_{t+\tau | t},
\]
where $\hat r_{t+\tau|t}$ is our prediction of the residual at time $t+\tau$
made at time $t$.
We form these predictions of future residual values using simple regression,
again on past and historical data.
Note that the baseline component of the prediction 
only depends on $t+\tau$, and not $t$, \ie,
the baseline value depends only on the time $t+\tau$ of the predicted quantity,
and not on the time $t$ at which the forecast is made.
The second term, our forecast of what the residual will be at time $t+\tau$,
does depend on $t$, the time of the forecast.

\subsection{Seasonal baseline}
The baseline is meant to capture the variation of the time series due to
time, typically, periodically repeating patterns.
A simple model for the baseline is a sum of $K$ sinusoids (\ie, Fourier terms),
\[
b_t = \beta_0 + \sum_{k = 1}^{K} \alpha_k \sin (2 \pi t / P_k) + \beta_k \cos (2 \pi t / P_k),
\] 
where $P_k$ are the periods.
Typically we would use as periods the fundamental period $P$ (\eg, one day, one year)
and those associated with the first few harmonics, \ie, $P/2$, $P/3$, \ldots.
We fit the coefficients $\alpha_k, \beta_k$, $k=1,\ldots, K$ using simple least squares
on historical data.

An example will illustrate the method.
Suppose the time period is fifteen minutes and we wish to model
diurnal (24 hour) and seasonal (annual) periodicities, with 4 harmonics each.
We choose
\[
P_1 = 96, \quad
P_2 = 48, \quad
P_3 = 24, \quad
P_4 = 12,
\]
as the periods for diurnal variation, and
\[
P_5 = 8766,  \quad
P_6 = 4383,  \quad
P_7 = 2101.5,  \quad
P_8 = 1095.75,
\]
for seasonal variation.
(One solar year is roughly 365 days and 6 hours, or 8766 periods of 15 minutes.)
This baseline model would have 17 parameters (including $\beta_0$, the constant).
If the time series is impacted by human or economic activity, 
we can also include weekly seasonality, or a weekend/holiday term.

Note that the value of the baseline model can be found for any time, 
past or future, once the baseline model coefficients are fixed, since 
it is then a fixed function of time.  We can, for example, evaluate the 
baseline load value, or renewable generation availablity,
at some specific time in the future.
In some applications of MPC, the baseline model is good enough to provide 
good performance.

\subsection{Auto-regressive residual forecasts}
Once we have the baseline forecast, we subtract it from our historical data to obtain
the sequence of historical residuals, $r_t = x_t - b_t$.
This sequence is sometimes referred to as the \emph{baseline adjusted} sequence.
(For example, with an annual baseline, $r_t$ is called the seasonally adjusted 
time series.)
Roughly speaking, $r_t$ contains the part of the sequence that is not explained 
by the periodicities.

To make forecasts of $r_{t+1}, \ldots, r_{t+T-1}$ at time $t$,
we use simple least squares regression based on the previous $M$ values,
$x_t, x_{t-1}, \ldots, x_{t-M+1}$.  Our model is
\[
 \hat r_{t+\tau|t} = \sum_{\tau' =0}^{M-1} \gamma_{\tau, \tau'} r_{t-\tau'},
\quad \tau = 1, \ldots, T-1,
\]
and we choose the $(T-1)\times M$ matrix of model 
parameters $\gamma_{\tau, \tau'}$ to minimize the mean square error on
the historical data.
These auto-regressive coefficients are 
readily interpretable: $\gamma_{\tau,\tau'}$ is the amount by
which $\hat r_{t+\tau|t}$ (our $\tau$-step-ahead prediction) depends on 
$r_{t-\tau'}$ (the value $\tau'$ steps in the past).

We can fit the coefficients associated with the forecast $\hat r_{t+\tau| t}$,
\ie, $\gamma_{\tau,\tau'}$ for $\tau' =0,\ldots, M-1$, separately for different
values of $\tau$.
Each of these is a separate least squares fit or regression, based on
historical data.
We note here a common error made in forecasting.
The bad method first
builds a `one-step-ahead' forecast, which gives $\hat r_{t+1|t}$.  Then, 
to forecast two steps ahead, the bad method iterates the one-step-ahead
forecast twice.  This method of iterating a one-step-ahead forecast
is more complicated,
and produces far worse forecasts, compared to the method
described above.

\subsection{Forecasting}
In summary, at time $t$ 
we predict the future values of the time series as
\[
 \hat x_{t+\tau|t} = b_{t+\tau} + \sum_{\tau'=0}^{M-1} \gamma_{\tau, \tau'} 
(x_{t-\tau'} - b_{t-\tau'}),
\quad \tau = 1, \ldots, T-1.
\]
This forecast depends on the baseline model coefficients, as well as
the residual auto-regressive coefficients.

There are many useful extensions of the basic forecasting method described above.
First, the baseline
could form predictions based not just on the time, but also on other observable
quantities, for example, weather.
In this case the baseline is not known ahead of time; it depends on the 
values of the other observable quantities at that time.
Another simple and useful extension is to project the forecast values onto
some set of values that the quantity is known to have.
For example if we are forecasting a renewable generator power
availability, which must lie between $0$ and some upper limit (the capacity)
$P^\mathrm{max}$, we clip our forecast value to lie within this range.
(That is, we replace the forecast value with $0$ if it is negative,
and $P^\mathrm{max}$ if it is larger than $P^\mathrm{max}$.)

While we have described the construction of the forecast as a two step process,
\ie, fitting a baseline, and then fitting an auto-regressive model for
the residuals, the two steps can in fact be done at the same time.  We simply
fit a predictor of $x_{t+\tau}$ for each $\tau$, using a set of regressors that
include current and previous values, the baseline basis functions, and indeed
any other regressors that might help, \eg, weather or futures contract prices. 
(That approach would give a predictor very 
close, but not equal, to the one described here.)
We have described the construction of the forecast as a two-step process because 
it is easy to interpret.

\subsection{Generating sample trajectories}
In the forecasting method described above, the goal is to come up with one
estimate of the future of the time series.
In this section we describe a simple method for generating a set of 
sample forecast trajectories
\[
x_{t+\tau|t}^{(k)}, \quad k=1, \ldots, K.
\]
These sample trajectories can be used for robust MPC, as 
described in \S\ref{s-robust-mpc-def}.
They are also useful as a sanity check on our forecasting method.
If the generated sample forecasts don't look right, it casts some 
doubt on our forecasting method.  If our forecasts look plausible, we gain
confidence in our forecast method.
The method we describe here works with \emph{any} forecasting method, 
including even the simplest ones, such as forecasting the value
as the current value, \ie, $\hat x_{t+\tau|t} = x_t$.

We let $e_t\in \reals^T$ denote the vector of forecast errors for
$x_{t+\tau|t}$, \ie,
\[
(e_t)_\tau = x_{t+\tau} - \hat x_{t+\tau|t}, \quad \tau = 0, \ldots, T-1.
\]
(For simplicity we index the vectors $e_t$ from $0$ to $T-1$.)
We collect these forecast error vectors over a historical data set,
and then fit these vectors with a Gaussian distribution $\mathcal N(\mu,\Sigma)$.
The simplest method uses the empirical mean and covariance of $e_t$ over
the historical data as $\mu$ and $\Sigma$, and in many cases, 
we can take $\mu=0$.
More sophisticated methods for choosing $\mu$ and $\Sigma$ involve 
adding a regularization term, or fitting a low-rank model.

To generate $K$ forecasts at time $t$, we sample $K$ vectors $e_t^{(k)}$,
$k=1, \ldots, K$ from $\mathcal N(\mu,\Sigma)$, and then form the sample 
forecasts as
\[
\hat x_{t+\tau|t}^{(k)} = 
\hat x_{t+\tau|t} + e^{(k)}_\tau, \quad \tau = 0, \ldots, T-1, \quad
k=1, \ldots, K.
\]
These samples are meant to be plausible guesses as what the next $T-1$ values 
of the time series might be.  Their average value is our forecast. We add 
to our forecast the simulated forecast errors that have the same mean and 
covariance of historically observed values.

\subsection{Wind farm example}
\label{s-mpc-ar-model}

We consider a time series of the power output
 of a wind farm, in MW,
(which depends on the available wind force)
on data by the National Renewable Energy  
Laboratory (NREL) for a site in West Texas.
The code can be seen in the Python notebook at
\url{https://github.com/cvxgrp/cvxpower/blob/master/examples/WindForecast.ipynb}.
Observations are taken every 5 minutes, from January 2010 to December 2012.
We use data from 2010 and 2011 to train the models, 
and data from 2012 for testing.
Our model has a baseline component that
uses 4 periodicities to model diurnal variation,
and the other 4 to model annual variation.
Our predictor uses an autoregressive predictor of the residuals
to the baseline (\ie, the seasonality-adjusted
series) of size $T = M = 288$, \ie, it forecasts
every 5 minutes power availability for the next 24 hours, 
using data from the past 24 hours.
Finally, 
since the power output lies between 0 and 16 MW (the minimum and maximum power 
of the turbine),
we project our forecast onto this interval.

\begin{figure*}
\begin{center}
\includegraphics[width=.9\textwidth]{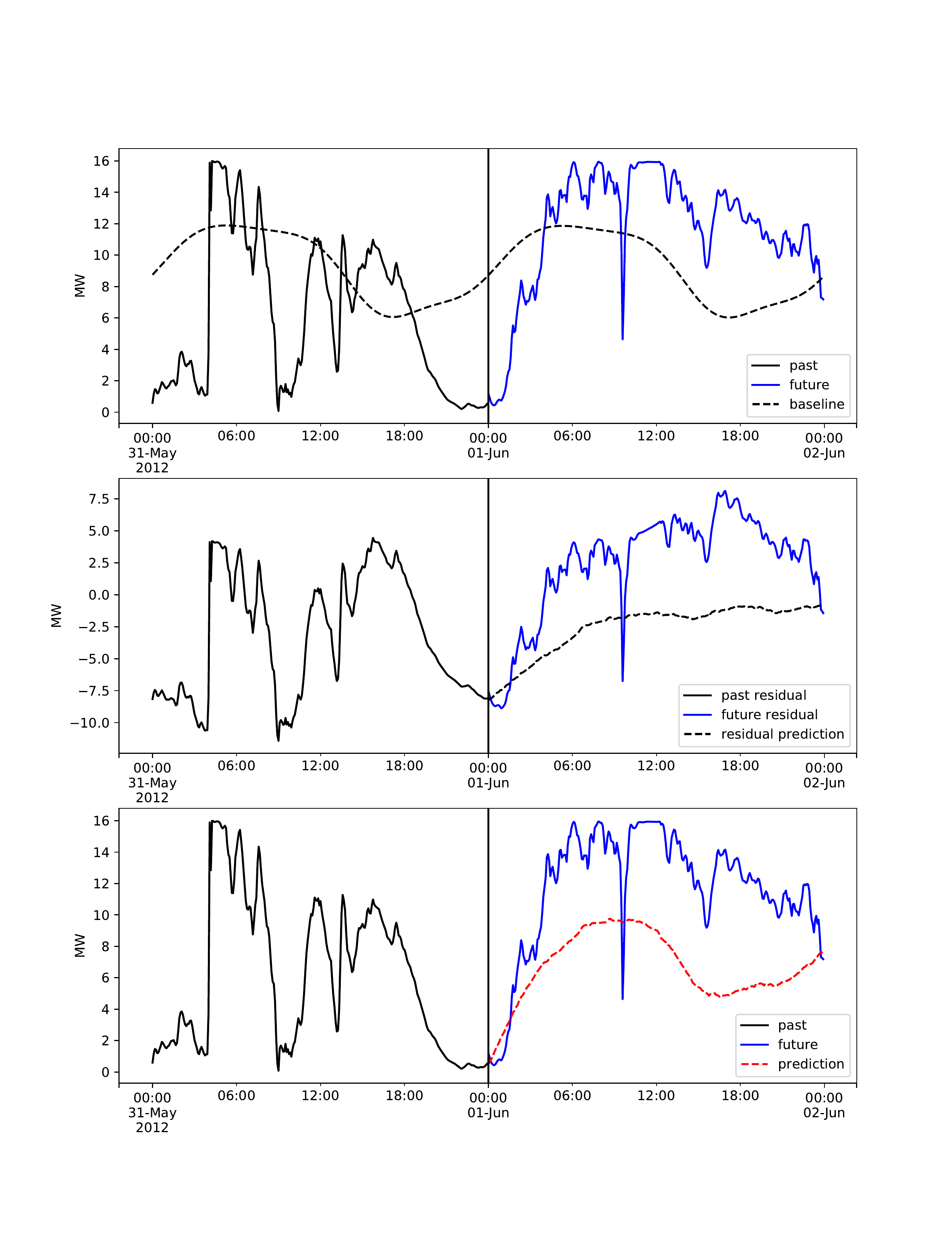}
\caption{
Wind farm power availability example.  Powers shown over two days,
with the vertical bar representing the time
the forecast is made.  \emph{Top.} True and baseline power.
\emph{Middle.} Residual and residual prediction.
\emph{Bottom.} True power and forecast.
}
\label{f-wind-forecast}
\end{center}
\end{figure*}

\begin{figure*}
\begin{center}
\includegraphics[width=.95\textwidth]{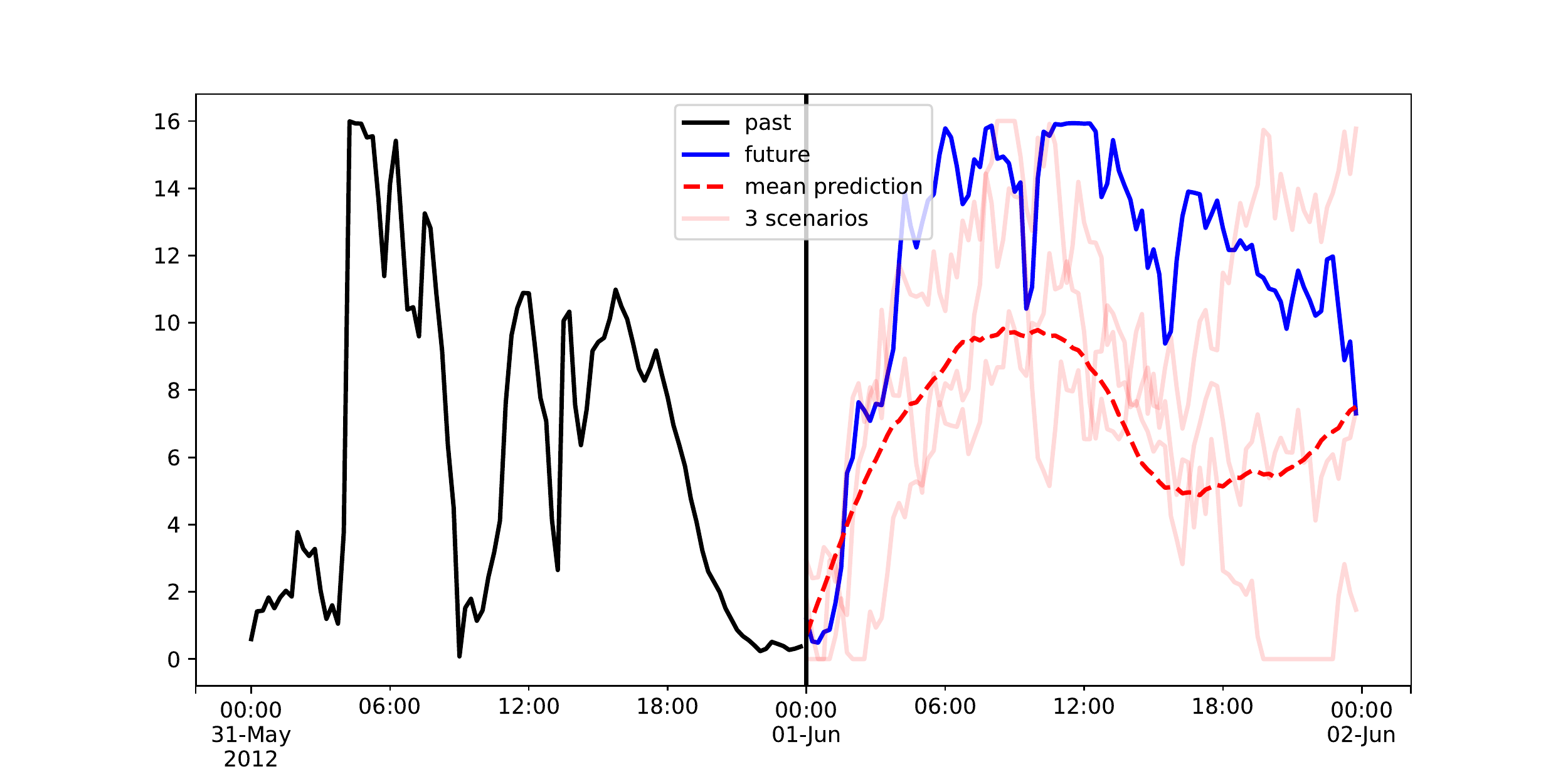}
\caption{Wind farm power availability example, forecast and 3 generated scenarios.}
\label{f-wind-many-scenarios}
\end{center}
\end{figure*}

Figure~\ref{f-wind-forecast} shows the result
of the forecast on a day in June 2012,
which is in the test set.
Figure~\ref{f-wind-many-scenarios} shows $K=3$ generated sample
trajectories, or scenarios, for the same day.  At least to the eye, they look
quite plausible.

\section{Appendix: Code example for \texttt{cvxpower}}\label{s-cvxpower}
We show here the Python source code to construct and optimize
the network of \S\ref{s-three-bus-example}.
We define objects for each load, generator,
transmission line, net, and then combine them
to formulate, and solve, the static optimal power flow
 problem. More examples (in the form of Python notebooks)
 can be seen in the ``examples'' folder of the software
 repository
 (at \url{https://github.com/cvxgrp/cvxpower}).

\begin{Verbatim}[samepage=true]
from cvxpower import * 

load1 = FixedLoad(power=50, name="load1")
load2 = FixedLoad(power=100, name="load2")

gen1 = Generator(power_max=1000, alpha=0.02, beta=30, name="gen1")
gen2 = Generator(power_max=100, alpha=0.2, beta=0, name="gen2")

line1 = TransmissionLine(power_max=50, name='line1')
line2 = TransmissionLine(power_max=10, name='line2')
line3 = TransmissionLine(power_max=50, name='line3')

net1 = Net([load1.terminals[0], gen1.terminals[0], 
            line1.terminals[0], line2.terminals[0]], name = 'net1')
net2 = Net([load2.terminals[0], line1.terminals[1], 
            line3.terminals[0]], name = 'net2')
net3 = Net([gen2.terminals[0], line2.terminals[1], 
            line3.terminals[1]], name = 'net3')
network = Group([load1, load2, gen1, gen2,
    line1, line2, line3],
    [net1, net2, net3])

network.init_problem()
network.optimize()
network.results.summary()
\end{Verbatim}
The output is: 
\begin{Verbatim}
Terminal                  Power
--------                  -----
load1[0]                  50.00
load2[0]                 100.00
gen1[0]                  -90.00
gen2[0]                  -60.00
line1[0]                  50.00
line1[1]                 -50.00
line2[0]                 -10.00
line2[1]                  10.00
line3[0]                 -50.00
line3[1]                  50.00

Net                       Price
---                       -----
net1                    33.6000
net2                   199.6002
net3                    24.0012

Device                  Payment
------                  -------
load1                   1680.00
load2                  19960.02
gen1                   -3024.00
gen2                   -1440.07
line1                  -8300.01
line2                    -95.99
line3                  -8779.95
\end{Verbatim}

\end{document}